%
%

\input amstex
\documentstyle{amsppt}
\pagewidth{6.4in}
\vsize8.5in
\parindent=6mm
\parskip=3pt
\baselineskip=14pt
\tolerance=10000
\hbadness=500
\NoRunningHeads
\loadbold
\topmatter
\title 
Homogeneous Fourier multipliers
of Marcinkiewicz type
\endtitle
\author Anthony Carbery and Andreas Seeger\endauthor
\thanks
The second author was supported in
part by a grant from the National Science Foundation (USA)
\endthanks
\address
Centre for mathematical analysis and its applications,
University of Sussex, Falmer, Brighton BN19QH, U. K. 
\endaddress
\address Department of mathematics, 
University of Wisconsin,  Madison, WI 53706, USA
\endaddress
\endtopmatter


\define\R{{\Bbb R}}

\define\ep{\epsilon}
\define\supp{{\text{\rm supp }}}

\define\inn#1#2{\langle#1,#2\rangle}

\define\la{\lambda}
\define\om{\omega}

\define\dlk{{|\lambda_{\kappa}(s)|ds}}
\define\k{\kappa}
\define\ka{\kappa}
\define\vphi{\varphi}
\define\tbeta{\widetilde \beta}

\define\fa{{\frak a}}

\redefine\L{{\Cal L}}

\head
{\bf 1. Introduction}
\endhead

Let $m\in L^\infty(\Bbb R^2)$ be homogeneous of degree zero. Then $m$ is 
almost everywhere determined by 
$h_\pm(\xi_1)=m(\xi_1,\pm 1)$.
For $k\in \Bbb Z$ let $I_k=[2^{-k-1},2^{-k}]\cup[-2^{-k},-2^{-k-1}]$ and 
let $h_+$ and $h_-$ satisfy the condition
$$
\sup_{k\in \Bbb Z} \,\Bigl(\int_{I_k}
\big |s h_\pm'(s)\big|^r \frac{ds}s\Bigr)^{1/r}\,<\,\infty.
\tag 1.1
$$
Rubio de Francia posed the question whether  a condition like
(1.1) is sufficient
to prove that $m$ is a Fourier multiplier of $L^p(\Bbb R^2)$, $1<p<\infty$. 
An application of the Marcinkiewicz multiplier theorem with
$L^2$-Sobolev hypotheses ({\it cf.} (1.3) and (1.5) below) and
interpolation arguments already show 
 that the answer is yes, provided $r>2$.
Recently, Duoandikoetxea and Moyua \cite{15} have shown that the same 
conclusion  can be reached if $r=2$.
On the other hand, since  characteristic functions of halfspaces
are Fourier multipliers of $L^p$, $1<p<\infty$,
 a simple averaging argument  shows that 
the condition 
$h'\in L^1$ implies $L^p$-boundedness for
$1<p<\infty$.
Our first theorem 
shows that the weaker assumption (1.1) with $r=1$ implies
boundedness in $L^p(\Bbb R^2)$, for
$1<p<\infty$.

\proclaim{ Theorem 1.1}
Suppose that $h_+$ and $h_-$ satisfy  the hypotheses
 of the Marcinkiewicz multiplier
theorem on the real line, that is 
$$\sup_{k\in \Bbb Z}\int_{I_k}|dh_\pm(s)|\,\le\, A
\tag 1.2
$$
for 
$I_k=[2^{-k-1},2^{-k}]\cup[-2^{-k},-2^{-k-1}]$.
Let $m\in L^\infty(\Bbb R^2)$ be homogeneous of degree zero, such that
for $\xi_1\in\Bbb R$,
$m(\xi_1,1)=h_+(\xi_1)$  and
$m(\xi_1,-1)=h_-(\xi_1)$.
Then $m$ is a Fourier multiplier of $L^p(\Bbb R^2)$,
 $1<p<\infty$, with norm $\le C\, A$.
\endproclaim

One can obtain
a stronger result for fixed  $p>1$ 
using the space $V^q$ of functions of 
bounded $q$-variation. Given an interval $I$ on the real line
a function $h$ belongs to $V^q(I)$ if for each partition
$\{x_0<x_1<\dots <x_N\}$ of $I$  the sum
$\sum_{\nu=1}^N|h(x_\nu)-h(x_{\nu-1})|^q$ is bounded and the upper bound 
of such sums is finite. We denote by $\|h\|_{V^q}^q$ the least upper bound.
Then the following result is an immediate consequence of Theorem 1.1  and the
interpolation argument in \cite{8}.

\proclaim{\bf Corollary 1.2}
Let  $m$,   $h_\pm$ and $I_k$ be as above
 and suppose that
$$\|h_\pm\|_\infty\,+\,\sup_k \|h_\pm\|_{V^q(I_k)}\,<\, \infty.$$
Then $m$ is a Fourier multiplier of $L^p(\Bbb R^2)$, if $|1/p-1/2|<1/2q$.
\endproclaim
 
A slightly weaker result can be formulated  in terms of Sobolev spaces. 
Let $\beta$ be  an even
$C^\infty$  function on the real line, supported in 
$(5/8,8/5)\cup (-8/5,-5/8)$ 
and positive in
$(1/\sqrt2,\sqrt2)\cup (-\sqrt2,-1/\sqrt2)$; we shall usually assume that 
$\sum_{k\in \Bbb Z}
\beta^2(2^k s)=1$ for $s\neq 0$.
Let 
$L^q_\alpha(\Bbb R^d)$ denote the standard
Sobolev space
with norm
$\|h\|_{L^q_\alpha}=\|\Cal F^{-1}[(1+|\xi|^2)^{\alpha/2}\widehat h]\|_q$.
Then $L^q_\alpha(\Bbb R)\subset V_q$ if $\alpha>1/q$ and therefore we obtain

\proclaim{Corollary 1.3}
Let $m\in L^\infty(\Bbb R^2)$   
be homogeneous of degree zero and  $h_\pm(\xi_1)=m(\xi_1,\pm 1)$. 
Suppose that  $q>1$ and that
$$\sup_{t\in\Bbb R_{+}}
\|\beta\, h_\pm(t\cdot)
\|_{L^q_\alpha(\Bbb R)}\,<\,\infty,\qquad
\alpha>\frac 1q.
\tag 1.3
$$
Then $m$ is a Fourier multiplier of $L^p(\Bbb R^2)$ if $|1/p-1/2|<1/2q$.
\endproclaim

We now compare these results with more standard
multiparameter versions of
the H\"ormander-Marcinkiewicz multiplier theorem.
In order to formulate them
let
$$\Cal D_j^\alpha g\,=\,\Cal F^{-1}[(1+|\xi_j|^2)^{\alpha/2}\Cal Fg]$$
and, for $1<q<\infty$, let 
$ {\Cal H}^q_\alpha(\Bbb R^{n})$
 be the {\it multiparameter} 
Sobolev space of all functions $g$, such that 
$$
\|g\|_{{\Cal H}^q_\alpha(\Bbb R^{n})}\,:=\,
\|\Cal D_1^\alpha\dots\Cal D_{n}^\alpha g\|_{L^q(\Bbb R^n)}\,<\,\infty.
$$
Let $\beta$ be  as above and 
denote by $\beta_{(i)}$ a copy of $\beta$ as a function of
the $\xi_i$-variable.
Then if $q\ge 2$ the condition
$$\sup_{t\in(\Bbb R_{+})^{d}}
\|\beta_{(1)}\otimes\dots\otimes\beta_{(d)}\, m(t_1\cdot,\dots,t_{d}\cdot)
\|_{\Cal H^q_\alpha(\Bbb R^{d})}\,<\,\infty,\qquad
\alpha>\frac 1q
\tag 1.4
$$
implies that $m$ is a Fourier multiplier of $L^p$ for 
$|1/p-1/2|< 1/q$. For $q=2$ the proof of this result  
is a variant of 
Stein's proof of the H\"ormander multiplier theorem (see \cite{25, ch.IV})  
and the general case follows by an interpolation argument as in \cite{9}.
If we apply this result to  homogeneous multipliers and
set
$$
m(\xi',\pm 1)\,=\, g_\pm(\xi'),\qquad\xi'\in\Bbb R^{d-1}
\tag 1.5
$$ 
we obtain by a straightforward computation
\proclaim{Corollary 1.4}
Suppose that $r\ge 2$,
$$\sup_{t\in(\Bbb R_{+})^{d-1}}
\bigl\|\Cal D_1^{2\gamma}\Cal D_2^{\gamma}\dots\Cal D_{d-1}^\gamma
\bigl[
\beta_{(1)}\otimes\dots\otimes\beta_{(d-1)}
 g_\pm(t_1\cdot,\dots,t_{d-1}\cdot)\bigr]
\bigr\|_{L^r(\Bbb R^{d-1})}\,<\,\infty,\qquad
\gamma>\frac 1r,
\tag 1.6
$$
and that the condition analogous  to $(1.6)$ holds for all permutations of
the $(s_1,\dots,s_{d-1})$-variables.
Let $m$ be homogeneous of degree zero and related to $g_\pm$ by {\rm(1.5)}.
Then $m$ is a Fourier multiplier of $L^p(\Bbb R^d)$ if $|1/p-1/2|<1/r$.
\endproclaim

In two dimensions Corollary 1.4 says that if 
$\alpha>1/q$, $q\ge 1$, and $\beta g_\pm(t\cdot)\in\Cal 
H^{2q}_\alpha(\Bbb R)$, 
uniformly in $t>0$, 
then
$m$ is a Fourier multiplier of $L^p$ if 
$|1/p-1/2|<1/2q$.  Corollary 1.3 is stronger since 
  a compactly supported function in
$H^{2q}_\alpha(\Bbb R)$ belongs to
$H^{q}_\alpha(\Bbb R)$.

We are now going to discuss
variants  of Theorem 1.1 in  higher dimensions. 
First if $g_\pm\in \Cal H^q_\alpha(\Bbb R^{d-1})$, $\alpha>1/q$
and if $g_{\pm}$ are compactly supported in  $[1/2,2]^{d-1}$  then  the 
homogeneous extension
$m$ is a Fourier multiplier of $L^p(\Bbb R^d)$ if $|1/p-1/2|<1/2q$.
In fact by a simple averaging argument one sees that the condition
$g_\pm\in \Cal H^1_{1+\epsilon}$ implies that $m$ is an $L^1$ multiplier
and the general case follows by interpolation. We remark that if
$\alpha<|2/p-1|$ the condition $g_\pm\in\Cal H^q_\alpha$ (any $q$)
does not imply that
$m$ is a Fourier multiplier of $L^p$.
Relevant counterexamples have been pointed out by L\'opez-Melero
\cite{22} and Christ \cite{7}.

Perhaps surprisingly, the situation in higher dimensions changes
 if one imposes
dilation invariant  conditions as in Theorem 1.1.
One might want to just replace hypothesis (1.2) by the hypotheses of the 
Marcinkiewicz multiplier theorem in $\Bbb R^{d-1}$
(\cite{25, p.108}). However this assumption is not sufficient
to deduce that $m$ is a Fourier multiplier of $L^p$  for any $p\neq 2$
(see \S3 for the counterexample involving the Kakeya set).
 However we do have

\proclaim{Theorem 1.5} 
 Let $m\in L^{\infty}(\R^d)$, $d\ge 2$, be homogeneous of degree zero and 
let $g_\pm$ be as in {\rm(1.5)}.
Suppose that $q\ge 2$, and
$$\sup_{t\in(\Bbb R_{+})^{d-1}}
\|\beta_{(1)}\otimes\dots\otimes\beta_{(d-1)}
\, g_\pm(t_1\cdot,\dots,t_{d-1}\cdot)
\|_{\Cal H^q_\alpha(\Bbb R^{d-1})}\,<\,\infty,\qquad
\alpha>\frac 1q.
\tag 1.7
$$
Then $m$ is a Fourier multiplier of $L^p(\Bbb R^d)$ if $|1/p-1/2|<1/2q$.
\endproclaim
Interpolating  Theorem 1.5 with Corollary 1.4 (with $p$ close to $1$) yields
\proclaim{Corollary 1.6} 
Let $m\in L^{\infty}(\R^d)$, $d\ge 2$, be homogeneous of degree zero and 
let $g_\pm$ be as in {\rm(1.5)}. Suppose that $1<p<4/3$ and
$$\sup_{t\in(\Bbb R_{+})^{d-1}}
\bigl\|\Cal D_1^{\alpha}\Cal D_2^{\gamma}\dots\Cal D_{d-1}^\gamma
\bigl[
\beta_{(1)}\otimes\dots\otimes\beta_{(d-1)}
 g_\pm(t_1\cdot,\dots,t_{d-1}\cdot)\bigr]
\bigr\|_{L^2(\Bbb R^{d-1})}\,<\,\infty,\qquad
\gamma>\frac 12,\quad \alpha>\frac 2p -1
$$
and that the analogous conditions obtained by permuting 
the $(s_1,\dots,s_{d-1})$-variables hold.
Then
$m$ is a Fourier multiplier of $L^p(\Bbb R^d)$.
\endproclaim
In particular if 
$\sup_{t\in(\Bbb R_{+})^{d-1}}
\|\beta_{(1)}\otimes\dots\otimes\beta_{(d-1)}
\, g_\pm(t_1\cdot,\dots,t_{d-1}\cdot)
\|_{\Cal H^2_\alpha(\Bbb R^{d-1})}\,<\,\infty$ and $1<p<4/3$  then $m$ 
is a Fourier multiplier of $L^p$ in $\alpha>\tfrac 2p-1$. This result is 
essentially sharp:
in \S3 we show that in order for
$\sup_{t\in(\Bbb R_{+})^{d-1}}
\|\beta_{(1)}\otimes\dots\otimes\beta_{(d-1)}
\, g_\pm(t_1\cdot,\dots,t_{d-1}\cdot)
\|_{\Cal H^q_\alpha(\Bbb R^{d-1})}\,<\,\infty$ to imply that $m$ is 
a Fourier multiplier of $L^p$  we must necessarily
have $\alpha\ge 2/p-3/2+1/q$ if $1<p<4/3$ and $\alpha>1/q$ if $4/3\le p\le 2$.

In order to prove more refined results  on $L^p(\Bbb R^d)$, $d\ge 3$, $p$
close  to $1$,
we shall use multiparameter Calder\'on-Zygmund theory.
It turns out that  it is useful (and easier) to first 
prove  a result for the multiparameter Hardy-space
 $H^p(\Bbb R^d)$, $0<p\le 1$. 
The Hardy space $H^p$ is defined in terms of square-functions invariant under
 the multiparameter family of dilations
$\delta_tx= (t_1x_1,\dots,t_dx_d)$, $t\in(\Bbb R_+)^d$.
Again we formulate  the multiplier result using localized
multiparameter Sobolev spaces 
invariant under 
multiparameter dilations. In order to include a sharp result also  for
$p<1$ we want to admit 
values of $q\le 1$ in (1.2). To make this possible the  definition of 
$\Cal H^q_\alpha$ has to be modified.
We may always assume that $\beta$  above 
is
such that $\sum_{r\in \Bbb Z}\beta^2(2^{-r}s)=1$ for 
$s\neq 0$. Let $\psi_r=\beta^2(2^{-r}\cdot)$ if $r\ge 1$ and
$\psi_0=1-\sum_{r>0}\psi_r$. 
For $n=(n_1,\dots,n_{d-1})$, $n_i\ge 0$, $i=1,\dots,d-1$ set
$\psi_n(\xi_1,\dots,\xi_{d-1})=\prod_{i=1}^{d-1}\psi_{n_i}(\xi_i)$.
The decomposition 
$$g=\sum_{n\in(\Bbb N_0)^{d-1}}
 \widehat \psi_n*g$$
is referred to as the inhomogeneous Littlewood-Paley decomposition of
$\Bbb R^{d-1}$.
Then
$$\|g\|_{\Cal H^q_\alpha(\Bbb R^{d-1})}
\approx \Bigl\|\Bigl (
\sum_{n\in(\Bbb N_0)^{d-1}}
2^{2(n_1+\dots+n_{d-1})
\alpha}|\widehat\psi_n*g|^2\Bigr)^{1/2}
\Bigr\|_{L^q(\Bbb R^{d-1})}
\tag 1.9
$$
for $1<q<\infty$, and for $q\le 1$ we define ${\Cal H^q_\alpha(\Bbb R^{d-1})}$
as the space of tempered distributions for which the quasinorm on
the right hand side of (1.9)  is finite. In this paper we shall always have
$\alpha>1/q$; in this case 
$\Cal H^q_\alpha$ is embedded in $L^\infty$.
This and other properties of the spaces $\Cal H^q_\alpha$
may be proved by obvious modifications of the one-parameter case; for the 
latter we refer to \cite {27}.

\proclaim{\bf Theorem 1.7}
Let $m\in L^\infty(\Bbb R^d)$ be homogeneous of degree zero and related to
$g_\pm$ as in  {\rm (1.5). }
Suppose that $0<r\le 1$ and
$$\sup_{t\in(\Bbb R_{+})^{d-1}}
\|\beta_{(1)}\otimes\dots\otimes\beta_{(d-1)}
 g_\pm(t_1\cdot,\dots,t_{d-1}\cdot)
\|_{\Cal H^r_\alpha(\Bbb R^{d-1})}\,<\,\infty,\qquad
\alpha>\frac 2r-1.
\tag 1.10
$$
Moreover if $d\ge 3$ suppose  that
$$\sup_{t\in(\Bbb R_{+})^{d-1}}
\bigl\|\Cal D_1^{2\gamma}\Cal D_2^{\gamma}\dots\Cal D_{d-2}^\gamma
\bigl[
\beta_{(1)}\otimes\dots\otimes\beta_{(d-1)}
 g_\pm(t_1\cdot,\dots,t_{d-1}\cdot)\bigr]
\bigr\|_{L^2(\Bbb R^{d-1})}\,<\,\infty,\qquad
\gamma>\frac 1r-\frac 12
\tag 1.11
$$
and that the analogous conditions obtained by permuting 
the $(s_1,\dots,s_{d-1})$-variables hold.
Then $m$ is a Fourier multiplier of the multiparameter Hardy space
$H^p(\R^d)$, $r\le p<\infty$.
\endproclaim

Note that in two dimensions Theorem 1.7 is a natural extension of Corollary 1.4
to $H^p$-spaces in product domains.
The examples in \S3 show that in
higher dimensions additional assumptions such as (1.11) are necessary.
When $d\ge 3$, Theorem 1.7 with $r=1$  serves as a
substitute for Theorem 1.1. 
Notice that if $r=1$ condition (1.10) involves mixed derivatives in $L^1$
of order $d-1+\epsilon$, and condition (1.11) involves derivatives 
in $L^2$ up to order $(d-1+\epsilon)/2$.  In comparison
the hypotheses in Corollaries 1.3 and 1.6 involve $L^2$ derivatives up 
to order $(d+\ep)/2$ if $p$ is close to $1$. 
As a consequence we obtain the following analogue of Corollary 1.4,
formulated in terms of the standard {\it oneparameter} Sobolev space
$L^q_\alpha$.

\proclaim{Corollary 1.8}
Let $m\in L^\infty(\Bbb R^d)$  
be homogeneous of degree zero and related to $g_\pm$ by {\rm(1.5)}.
Suppose that  $q>1$ and that
$$\sup_{t\in(\Bbb R_{+})^{d-1}}
\|\beta_{(1)}\otimes\dots\otimes\beta_{(d-1)}
\, g_\pm(t_1\cdot,\dots,t_{d-1}\cdot)
\|_{L^q_\alpha(\Bbb R^{d-1})}\,<\,\infty,\qquad
\alpha>\frac {d-1}q.
$$
Then $m$ is a Fourier multiplier of $L^p(\Bbb R^2)$ if $|1/p-1/2|<1/2q$.
\endproclaim

The counterexamples in \cite{22}, \cite{7} show
that the statement of the Corollary is false in the range
$|1/p-1/2|>1/2q$. However in view of Theorems 1.5 and 1.7
one expects the following sharper result. Namely suppose that for some
$q\in(1,2]$
$$\sup_{t\in(\Bbb R_{+})^{d-1}}
\|\beta_{(1)}\otimes\dots\otimes\beta_{(d-1)}
\, g_\pm(t_1\cdot,\dots,t_{d-1}\cdot)
\|_{\Cal H^q_\alpha(\Bbb R^{d-1})}\,<\,\infty,\qquad
\alpha>\frac 1q,
\tag 1.12
$$
and in dimension $d\ge 3$ suppose that 
$$\sup_{t\in(\Bbb R_{+})^{d-1}}
\bigl\|\Cal D_1^{\alpha}\Cal D_2^{\gamma}\dots\Cal D_{d-2}^\gamma
\bigl[
\beta_{(1)}\otimes\dots\otimes\beta_{(d-1)}
 g_\pm(t_1\cdot,\dots,t_{d-1}\cdot)\bigr]
\bigr\|_{L^2(\Bbb R^{d-1})}\,<\,\infty,\qquad
\gamma>\frac 12,\quad \alpha>\frac 1q
\tag 1.13
$$
as well as the analogous conditions obtained by permuting 
the $(s_1,\dots,s_{d-1})$-variables.
Then $m$ should be  a Fourier 
multiplier of $L^p(\Bbb R^d)$ if $|1/p-1/2|<1/2q$.
In order to prove this one is tempted
to use analytic interpolation and interpolate between the
$L^{p_0}$-estimate of  Theorem  1.7, for $p_0$ close to $1$, and the
$L^{4/3}$-estimate of Theorem 1.5.
One would have to find the intermediate spaces for
intersections of $L^2$ and $L^q$ Sobolev spaces.
However the intersection of the intermediate 
spaces does not need to be contained 
in the intermediate space of the intersections
(for related counterexamples see \cite{26}). It is actually possible to prove 
the result for $|1/p-1/2|<1/2q$ (assuming (1.12), (1.13)) by
another approach. One has to use a general 
theorem for analytic families of 
operators acting on various kinds of atoms the proof of which relies
heavily on multiparameter Calder\'on-Zygmund theory. 
We do not include the technical proof here but refer the reader to
\cite{5}.

The paper is organized as follows: 
In \S2 we prove Theorem 1.1 using weighted norm inequalities and
variants of the maximal operator with respect to lacunary directions.
Examples demonstrating the sharpness of our results in higher 
dimensions  are discussed in \S3. 
The proof of Theorem 1.5 is in \S4; it relies on weighted norm inequalities 
which involve variants of the Kakeya maximal function. 
In \S5 we prove the Hardy space estimates of Theorem 1.7.

As a convention we shall refer to the quasi-norms in $H^p$ and 
$\Cal H^p_\alpha$ as ``norms'' although for $p<1$ these spaces are not 
normed spaces. 
By $M_p$, $1\le p\le \infty$, we denote 
the standard space of
Fourier multipliers of $L^p$. 
It will always be assumed that the even function $\beta\in C^\infty_0$ 
defined above  satisfies
$\sum_{\k\in\Bbb Z}[\beta(2^\k s)]^2=1$ for $s\neq 0$.
If $\fa\in\{1,\dots,d\}$ and $k$, $\tilde k$ in 
$\Bbb R^\fa$ then we shall use the notation
$k\le\tilde k$ if $k_i\le\tilde k_i$ for all 
$i\in\fa$. Similarly define $k\le \tilde k$ etc.
$C$ will always be an abstract constant which may assume
different values in different lines.

\bigskip

\head
{\bf 2. $\boldkey{L}^{\boldkey p}$-estimates in the plane}
\endhead

In the  proof of Theorem 1.1 there is no loss of generality 
in assuming that
$m$ is supported in the quadrant where
$\xi_1>0$ and $\xi_2>0$. By a limiting
argument as in Stein's book
\cite{25, p.109},  it suffices to prove the theorem
under the formally stronger assumption
$$\|h\|_\infty+
\sup_{k\in\Bbb Z}\int_{I_k}|h'(s)|\,ds\,\le\, A.
$$
Let $\beta$ be the smooth bump function defined in the introduction
(supported in $\pm[5/8,8/5]$). Let $\gamma \in C^\infty(\Bbb R^2\setminus
\{0\})$  be homogeneous of degree $1$ such that $\gamma(\xi)=1$ if
$|\xi_1/\xi_2|\in[25/64,64/25]$ (in particular on the 
support of $\beta\otimes \beta$) and such that 
 $\gamma(\xi)=0$ if
$|\xi_1/\xi_2|\notin(1/4,4)$.
Set
$$h_\ka(\xi)=\gamma(\xi)h(2^{-\kappa}\xi_1/\xi_2).
$$
Then we may split
$$m\,=\,\sum_{k\in\Bbb Z^2}[\beta \otimes\beta
m_k](2^{k_1}\cdot, 2^{k_2}\cdot)$$
where
$$
\align
m_k(\xi)\,&=\,\beta(\xi_1)\beta(\xi_2)h_{k_1-k_2}(\xi_1/\xi_2)
\\
&=\,\beta(\xi_1)\beta(\xi_2)
\int_0^{\xi_1/\xi_2}h_{k_1-k_2}'(s)\,ds.
\tag 2.1
\endalign
$$
and  $h_\k$ is supported in $(1/4,4)$, for all $\k\in\Bbb Z$.
Also set
$$
\widehat{T_kf}(\xi)\,=\,
[\beta \otimes\beta\,
m_k]
(2^{k_1}\xi_1, 2^{k_2}\xi_2)]\widehat f(\xi).
$$
Then by standard multiparameter Littlewood-Paley theory and duality, 
to establish 
Theorem 1.1 for $p\in [2,p_0)$, $p_0<\infty$,
it suffices to obtain an inequality 
$$
\int|T_k f|^2\omega\,\le\,C\,A^2\int|f|^2\frak M\omega
\tag 2.2
$$ for a certain operator $\omega\mapsto \frak M\omega$ which is bounded on
$L^q(\Bbb R^2)$ for $(p_0/2)'<q\le\infty$.
By our assumption on $h$,
$$
\sup_{\k\in \Bbb Z}\int |h_\k'(s)|ds\,\le C A.
\tag 2.3
$$
We denote  by $\L_k$ the standard Littlewood-Paley operator, such that
$$\widehat {\L_k f}(\xi)\,=\,\beta(2^{k_1}\xi_1)\beta(2^{k_2}\xi_2)
\widehat f(\xi)$$
and define the operator 
$S_{ks}$ by
$$\widehat{S_{\k s}f}(\xi)\,=\,
\cases
\widehat f(\xi),\qquad &\text{ if }2^{\k}\xi_1/\xi_2>s,\,
\xi_1\ge 0,\,\xi_2\ge 0
\\
0,\qquad &\text{ otherwise}
\endcases .
$$
Then from (2.1) we  see that
$$
T_k f(x)\,=\,\int_{1/8}^8
\L_k S_{k_1-k_2,s}f(x) h_{k_1-k_2}'(s) ds
$$
Then,
if $\omega\ge 0$ is a weight we apply the Cauchy-Schwarz inequality
to obtain
$$
\int|T_k f(x)|^2\omega(x)\,dx\,\le\,
C\,A\,\iint |\L_k S_{k_1-k_2,s}f(x)|^2 |h_{k_1-k_2}'(s)| ds
\,\omega(x)dx.
\tag 2.4
$$
Let $M_{(1)}$, $M_{(2)}$ be the Hardy-Littlewood maximal
functions with respect to the coordinate directions and let
$M_{\k,s}$ be the 
Hardy-Littlewood maximal
function with respect to the  direction perpendicular to
$\{\xi;\,2^{\k}\xi_1/\xi_2=s\}$, {\it i.e.} in the direction
$(1, -2^{-\k}s)$. Then using 
 weighted norm inequalities for singular integral operators
due to C\'ordoba and Fefferman (\cite{13}, see also \cite{18})
we see that the expression on the right hand side of (2.4) is dominated by
$$
C_\alpha A\,
\int|f(x)|^2 M_{(1)}M_{(2)}
\Bigl[\int_{1/8}^8 (M_{k_1-k_2,s}\omega^\alpha)^{1/\alpha}
|h_{k_1-k_2}'(s)|ds\Bigr](x)\, dx
$$
where $\alpha>1$.
Now the proof of (2.2)
is completed by the following 

\proclaim{Proposition 2.1}
Let, for $\alpha\ge 1$,
$$\frak M_\alpha\omega(x)\,=\,\sup_{\k\in\Bbb Z}
\int_I(M_{\k,s}\omega^\alpha)^{1/\alpha}(x)|\la_\k(s)|ds
$$
where $I=[1/8,8]$ and
$$\sup_{\k\in \Bbb Z} \int_I|\la_\k(s)|ds\,\le\,B<\infty.$$
Then $\frak M_\alpha$ is bounded on $L^p(\Bbb R^2)$, $\alpha<p<\infty$,
with norm $\le C_{p,\alpha}B$.
\endproclaim

\demo{\bf Proof}
Since 
$$\frak M_\alpha(\omega)\,\le\,B^{1-\frac{1}{\alpha}}
[\frak M_1(\omega^\alpha)]^{\frac 1\alpha}
$$
it suffices to prove that $\frak M_1$ is bounded on $L^p$, $1<p<\infty$
with norm $C_p B$. If $\alpha>1$ then
$\frak M_\alpha$ will be bounded  on $L^p$, $p>\alpha$, with norm
$\le C_{p/\alpha}^{1/\alpha}B$.

We follow arguments by Nagel, Stein and Wainger \cite{23} as modified by
Christ (see \cite{2}).
Let $\vphi:\Bbb R\to\Bbb R$ be smooth, even, nonnegative,
with  $\vphi(0)>0$
such that $\widehat\vphi$ has compact support in
$[-1/20,1/20]$. Let
$$
\psi(\xi_1,\xi_2)\,=\,\widehat\vphi(\xi_1+\xi_2)
$$
and define for $\k\in \Bbb Z$
$$\widehat{P^l_{\k s}\omega}(\xi)\,=\,
\psi(2^l \xi_1,s 2^{l-\k}\xi_2)\widehat\omega(\xi).
$$
It suffices to show that for $1<p<\infty$, $N$ being an arbitrary 
positive integer
$$
\biggl\|\sup_{-N\le\k\le N}\int_I
\sup_{l\in \Bbb Z}|P^l_{\k s}\omega|\,\dlk
\biggr\|_p\,\le\, C_p B\|\om\|_p
\tag 2.5
$$
where $C_p$ is independent of $N$. Then an application of the monotone 
convergence theorem allows to pass to the limit.
We note that for fixed $\k$, and $1<p\le\infty$
$$
\align
\biggl\|\int_I
\sup_{l\in \Bbb Z}|P^l_{\k s}\omega|\,\dlk
\biggr\|_p\,
\,&\le\,C\int_I\bigl\|\sup_{l\in\Bbb Z}|P^l_{\k s}\omega|\bigr\|_p\dlk
\\\,&\le\,C_p\,B\|\om\|_p
\tag 2.6
\endalign
$$
by the $L^p$ estimate for the one-dimensional
 Hardy-Littlewood maximal function $M_{k,s}$.
This means that we know a priori that the left hand side of (2.5) is bounded by
$B\,C_p(N)\|f\|_p$ (with $C_p(N)\le C_p'N$)
 and it remains to be shown that $C_p(N)$ can be chosen independently of
$N$. In what follows we define $C_p(N)$ to be the best constant in (2.5).

We first consider the case $2\le p<\infty$. 
Since the $L^\infty$-estimate is trivial  it suffices to prove the $L^2$
inequality.
We smoothly split $\psi$ into two parts, $\psi=\psi^0+\psi^1$ with
$\psi^1$ supported in the unit ball 
and $\psi^0$ supported in the cone
$\{\xi;\,|\xi_1+\xi_2|/|\xi|\le 1/2\}$.
We correspondingly define the operators $P^{l,0}_{\k s}$ and
$P^{l,1}_{\k s}$. 
Note that  there is the
 pointwise inequality
$$
|P^{l,1}_{\k s}\om(x)|\,\le\,C M_{(1)}M_{(2)}\om(x)
\tag 2.7
$$
which implies
$$
\biggl\|\sup_\k\int_I
\sup_{l\in \Bbb Z}|P^{l,1}_{\k s}\omega|\,\dlk
\biggr\|_p\,\le\, C_p B\|\om\|_p,\qquad 1<p\le \infty.
\tag 2.8
$$

Concerning $P^{l,0}_{\k s}$ 
we have
$$
|P^{l,0}_{\k s}\om(x)|\,\le\,C\bigl[ M_{(1)}M_{(2)}\om(x)+M_{\k,s}\om(x)\bigr]
$$
and therefore 
$$
\biggl\|\int_I
\sup_{l\in \Bbb Z}|P^{l,0}_{\k s}\omega|\,\dlk\biggr\|_p\,
\,\le\,C_{p,0}\,B\|\om\|_p
\tag 2.9
$$
for $1<p\le\infty$. Note that
$$\psi^0(2^l\cdot,s 2^{l-\k}\cdot)\,\widehat{ \om}\,=\,
\chi(\cdot,2^{-\k}\cdot)\psi^0(2^l\cdot,s 2^{l-\k}\cdot)\,\widehat{ \om}
$$
where $\chi$ is smooth, homogeneous of degree zero, identically $1$ on
$\{\xi;\,|\xi_1+\xi_2|/|\xi|\le 4\}$
and zero on $\{\xi;\,|\xi_1+\xi_2|/|\xi|\ge 8\}$.
Define the standard angular Littlewood-Paley operator $R_\k$ by
$$\widehat{R_\k f}(\xi)\,=\,
\chi(\xi_1,2^{-\k}\xi_2)\widehat f(\xi).$$
Then
$$
P^{l,0}_{\k s}\om\,=\, P^{l,0}_{\k s}R_\k \om
\tag 2.10
$$
and, as a consequence of multiparameter Littlewood-Paley theory and the
Marcinkiewicz multiplier theorem,
$$
\Bigl\|\bigl(\sum_\k|R_\k f|^2\bigr)^{1/2}\Bigr\|_p\,\le\,C\|f\|_p,
\qquad 1<p<\infty.
\tag 2.11
$$
Now by (2.10)
$$
\sup_{\k\in \Bbb Z}\Bigr|
\int_I\sup_{l\in\Bbb Z}|P^{l,0}_{\k s}\om|\dlk\Bigl|
\,\le\,\Bigl(
\sum_\k\Bigl[\int_I
\sup_{l\in\Bbb Z}|P^{l,0}_{\k s}R_\k\om|\dlk\Bigl]^2\Bigr)^{1/2}
\tag 2.12
$$
and using (2.9) we see that the
 square of the $L^2$-norm of the right hand side  equals
$$
\align
&\sum_\k\Bigl\|
\int_I\sup_{l\in\Bbb Z}|P^{l,0}_{\k s}R_\k\om|\dlk\Bigr\|_2^2
\\ \le\,
&\sum_\k\Bigl[\int_I\|\sup_{l\in\Bbb Z}|P^{l,0}_{\k s}R_\k\om\|_2\dlk\Bigr]^2
\\ \le\,&C B^2\sum_\k\|R_k \om\|^2_2\,\le\,C'B^2\,\|\om\|^2_2.
\endalign
$$
We have proved
$$
\Bigl\|\sup_{\k\in \Bbb Z}\Bigl[
\int_I\sup_{l\in\Bbb Z}|P^{l,0}_{\k s}\om|\dlk\Bigr]\Bigr\|_p
\,\le\,C B\|\om\|_p,\qquad 2\le p\le\infty.
\tag 2.13
$$
By (2.8) and (2.13) we see that
$$ C_p(N)\le C,\qquad 2\le p\le \infty.$$

We now assume $1<p<2$
and 
 begin with the observation that for any sequence $\{\omega_k\}$ of weights
we have
$$
\Bigl\|\Bigl(\sum_\k\Big|\int_I
\sup_{l\in \Bbb Z}|P^{l,0}_{\k s}\omega_\k|\,\dlk\Big|^p\Bigr)^{1/p}\Bigr\|_p\,
\le\, C_{p,0} B\bigl\|\bigl(\sum_\k|\om_\k|^p\bigr)^{1/p}\bigr\|_p,\qquad
1<p<\infty.
\tag  2.14
$$
This is immediate from (2.9).
Next positivity of $P^l_{ \k s}$ implies that
$$
\align
&\Bigl\|\sup_\k\int_I\sup_l|P^l_{\k s}\om_\k|\,\dlk\Bigr\|_p
\\
\le\,&
\Bigl\|\sup_\k\int_I\sup_l \big |P^l_{\k s}[\sup_\mu|\om_\mu|]\big|\,
\dlk\Bigr\|_p
\\
\,\le\,&B\,C_p(N)\bigl\|\sup_\k|\om_\k|\bigr\|_p
\tag 2.15
\endalign
$$
by the definition of
$C_p(N)$.
From (2.15) and (2.8) it follows that for  $1<p\le \infty$
$$
\Bigl\|\sup_{\k\in\Bbb Z}
\int_I\sup_{l\in\Bbb Z}|P^{l,0}_{\k s}\om_\k|\dlk\Bigr\|_p\,\le\,
C_{p,1} B\,C_p(N)\|\sup_\k [|\om_\k|]\|_p.
\tag 2.16
$$
Now if we interpolate (2.14) with (2.16) we obtain
for $p\le q\le\infty$
$$
\Bigl\|\Bigl(\sum_{\k\in\Bbb Z}\bigl|
\int_I\sup_{l\in\Bbb Z}|P^{l,0}_{\k s}\om_\k|\dlk\bigr|^q\Bigr)^{1/q}
\Bigr\|_p\,\le\,
C_{p,2} B\,C_p(N)^{1-p/q}\Bigl\|\Bigl(\sum_\k|\om_\k|^q\Bigr)^{1/q}
\Bigr\|_p.
\tag 2.17
$$
Using (2.12), (2.17) and (2.11) we obtain for $1 <p\le 2$
$$
\align
&\Bigl\|\sup_{\k\in \Bbb Z}\Bigr|
\int_I\sup_{l\in\Bbb Z}|P^{l,0}_{\k s}\om|\dlk\Bigl|\Bigr\|_p
\\ \le\,&\Bigl\|\Bigl(
\sum_\k\bigl[
\sup_{l\in\Bbb Z}|P^{l,0}_{\k s}R_\k\om|\dlk\bigl]^2\Bigr)^{1/2}\Bigr\|_p
\\
\le\,&C_{p,2} B\,C_p(N)^{1-p/2}
\Bigl\|\Bigl(\sum_\k|R_\k \om|^2\Bigr)^{1/2}\Bigr\|_p
\\
\le\,&C_{p,3} B\,C_p(N)^{1-p/2}
\|\om\|_p
\tag 2.18
\endalign
$$
Finally it follows from (2.8) and (2.18) that
$$
\,C_p(N)\,\le\, \bigl[C_p'+
C_{p,3} \,C_p(N)^{1-p/2}\bigr]
$$
which implies  that $C_p(N)$ is bounded 
by a constant depending only on $p$ but not on $N$. This finishes 
the proof of the proposition.\qed
\enddemo

\bigskip

\head
{\bf 3. Examples in higher dimensions}
\endhead
We show in this section  that Theorem 1.1 and Corollary 1.4 have no 
immediate analogue in terms of localized multiparameter Sobolev spaces
in higher dimensions. Our examples imply the sharpness of Theorems
1.5 and 1.7.

Let $L^p(L^2)$ be the space of functions $f$ on 
$\Bbb R^d=\Bbb R^{d_1}\oplus\Bbb R^{d_2}$
such that
$$
\|f\|_{L^p(L^2)}\,=\,\Bigl(\int\Bigl[\int |f(x',x'')|^2 dx''\Bigr]^{p/2}
dx'\Bigr)^{1/p}\,<\,\infty.
$$
For a bounded function $m$ on $\Bbb R^d$ we denote  by
$\|m\|_{M_p} $ the operator norm of the convolution operator $T$ 
defined by
$\widehat{T f}=m\widehat f$ and  by $\|m\|_{M_{p2}}$
the norm of $T$  as a bounded operator on $L^p(L^2)$.
 By a theorem of Herz and Rivi\`ere
\cite{19}
$$\|m\|_{M_{p,2}}\,\le\,C\|m\|_{M_p}
\tag 3.1
$$
for $1\le p\le\infty$.
We shall use the following 

\proclaim{Lemma 3.1}
Let $\{m_\k\}$ be a sequence  of bounded functions in $\Bbb R^{d_1}$.
Let $\chi\in C^\infty(\Bbb R^{d_2})$ 
be supported in $\{\xi'';1/2\le|\xi''|\le 2\}$ and equal to $1$ if
$1-\epsilon \le|\xi''|\le 1+\epsilon$ for some $\epsilon >0$.
Let
$$
m(\xi',\xi'')\,=\, \sum_\k\chi(2^{-6\k}\xi'')m_\k(\xi')
$$
and define $T_\ka$ by $\widehat{T_\ka f}(\xi)=m_\ka(\xi)\widehat f(\xi)$.
Then for
$1<p<\infty$
we have the inequality
$$
\Bigl\|\Bigl(\sum_\k|T_\k f_\k|^2\Bigr)^{1/2}\Bigr\|_{L^p(\Bbb R^{d_1})}
\,\le\,C_p
\|m\|_{M_{p2}(\Bbb R^d)}
\Bigl\|\Bigl(\sum_\k|f_\k|^2\Bigr)^{1/2}\Bigr\|_{L^p(\Bbb R^{d_1})}
$$
\endproclaim

\demo{\bf Proof}
Let  $\beta_0\in C^\infty$ be supported in 
$\{\xi'':1-\epsilon\le|\xi''|\le 1+\epsilon\}$ 
such that 
$$\|\beta_0\|_{L^2(\Bbb R^{d_2})}=1.$$
Let
$$g_\k(x',x'')\,=\, 2^{3\k}
 f_\k(x')\Cal F^{-1}_{\Bbb R^{d_2}}[\beta_0(2^{-6\k}\cdot)](x'')
$$
then by an application of  Plancherel's theorem  in the second variable
it follows that
$$
\Bigl\|\Bigl(\sum_\k|g_\k|^2\Bigr)^{1/2}\Bigr\|_{L^p(\Bbb R^{d})}\,=\,
(2\pi)^{-d_1/2}
\Bigl\|\Bigl(\sum_\k|f_\k|^2\Bigr)^{1/2}\Bigr\|_{L^p(\Bbb R^{d_1})}.
\tag 3.2
$$
Next let $L_\k $ denote convolution in $\Bbb R^{d_2}$ with
$\Cal F_{\Bbb R^{d_2}}^{-1}[\beta_0(2^{-\k}\cdot)]$.
By Littlewood-Paley theory we have for $1<p<\infty$
$$\Bigl\|\sum_\k L_\k g_\k\Bigr\|_{L^p(\Bbb R^d)}\,\le\, C_p
\Bigl\|\Bigl(\sum_\k|g_\k|^2\Bigr)^{1/2}\Bigr\|_{L^p(\Bbb R^{d})}.
\tag 3.3
$$
Now
$$
\align
\Bigl\|\Bigl(\sum_\k|T_\k f_\k|^2\Bigr)^{1/2}\Bigr\|_{L^p(\Bbb R^{d_1})}
\,&=\,
\Bigl\|\Bigl(\sum_\k
\int|\beta_0(2^{-6\k}\xi'')2^{-3\k}T_\k f_\k|^2d\xi''\Bigr)^{1/2} 
\Bigr\|_{L^p(\Bbb R^{d_1})}
\\
\,&=\,(2\pi)^{d_1/2} 
\Bigl\|\Bigl(\sum_\k\big|\Cal F_{\Bbb R^d}^{-1}[m_\k 
 \Cal F_{\Bbb R^d}
[ L_\k L_\k g_\k]]\big|^2\Bigr)^{1/2}
\Bigr\|_{L^p(L^2)}
\\
\,&=\,(2\pi)^{d_1/2} 
\Bigl\|\Bigl(\sum_\k\big|L_\k\Cal F_{\Bbb R^d}^{-1}\bigl[m  \Cal F_{\Bbb R^d}
[ \sum_j L_j g_j]\bigr]\big|^2\Bigl)^{1/2}\Bigr\|_{L^p(L^2)}
\endalign
$$
where the last identity holds in view of the support properties of $\beta_0$.
By Littlewood-Paley theory
$$
\align
\Bigl\|\Bigl(\sum_\k&\big|L_\k\Cal F_{\Bbb R^d}^{-1}\bigl[m  \Cal F_{\Bbb R^d}
[\sum_j L_j g_j]\bigr]\big|^2\Bigl)^{1/2}\Bigr\|_{L^p(L^2)}
\\
\,&\le\,
C_p \Bigl\|\Cal F_{\Bbb R^d}^{-1}\bigl[m  \Cal F_{\Bbb R^d}
[\sum_j L_j g_j]\bigr]\Bigr\|_{L^p(L^2)}
\\
\,&\le\,C_p\,\|m\|_{M_{p2}}
\Bigl\|\sum_j L_j g_j\Bigr\|_{L^p(L^2)}
\\
\,&\le\,C_p'\,\|m\|_{M_{p2}}
\Bigl\|\Bigl(\sum_j\bigl|  g_j\bigr|^2\Bigr)^{1/2}\Bigr\|_{L^p(L^2)}
\\
\,&=\,
C_p'\,\|m\|_{M_{p2}} (2\pi)^{-d_1/2}
\Bigl\|\Bigl(\sum_j|f_j|^2\Bigr)^{1/2}\Bigr\|_{L^p(\Bbb R^{d_1})}.\qed
\endalign
$$
\enddemo

We now show that the restriction $q\ge 2$ (corresponding to $4/3\le p\le 4$)
in   Theorem 1.5 is necessary.
 In what follows we denote by $L^p(L^2)$ the space of 
functions in $\Bbb R^3$ with
$$
\|f\|_{L^p(L^2)}\,=\,\Bigl(\iint\Bigl[\int |f(x_1,x_2,x_3)|^2 dx_2\Bigr]^{p/2}
dx_1 dx_3\Bigr)^{1/p}\,<\,\infty
$$
and correspondingly define 
$M_{p2}$.

Fix $N\gg 0$ and let
$$g_{N}(s_1,s_2)\,=\,
\sum_{\k=2}^{N}  \eta(N(s_1-\alpha_\k))
\widetilde\chi(2^{-6\k}s_2);
\tag 3.4
$$
where
$$
\alpha_\k\,=\,1+\frac{2N}{\k},
\tag 3.5
$$
and $\eta\in C^\infty$ is nonnegative, equal to $1$
in $[-1/4,1/4]$ and supported in $[-1/2,1/2]$. 
Similarly $\widetilde\chi$ is as in Lemma
3.1, supported in $\pm (1/2,2)$ and equal to $1$ in 
$\pm (1/\sqrt 2,\sqrt 2)$.
Then
$$
\sup_{s_1,s_2>0}\|\beta_{(1)}\otimes\beta_{(2)}g_{N}(s_1\cdot,s_2\cdot)
\|_{\Cal H^q_\alpha(\Bbb R^2)}\,\le\,
CN^{\alpha-1/q}.
\tag 3.6
$$

\proclaim{Lemma 3.2} 
Let $m_{(N)}$  be the homogeneous extension of $g_{N}$ 
defined in $(3.4)$, $(3.5)$.
There is a positive constant $c$ such that
$$\|m_{(N)}\|_{M_{p,2}}
\,\ge\,
\cases
c N^{1/2-2/p},\qquad&4<p<\infty
\\
c (\log N)^{1/4},\qquad &p=4.
\endcases
\tag 3.7
$$
\endproclaim

A comparison of (3.6) and (3.7) shows that in the case $p>4$ the condition
$$\sup_{s_1,s_2>0}\|\beta_{(1)}\otimes\beta_{(2)}g_\pm(s_1\cdot,s_2\cdot)
\|_{\Cal H^q_\alpha(\Bbb R^2)}\,<\,\infty
$$
does not imply $m\in M_p$ 
for the homogeneous extension $m$ if 
$\alpha<1/2+1/q-2/p$ (it does not even imply $m\in M_{p2}$). 
Similar statements follow by duality
for $1<p<4/3$. This yields the sharpness of Theorem 1.5.
By interpolation 
 an improvement of the $H^p$ estimates would lead 
to an improvement of the $L^p$ estimates and this implies
the sharpness of Theorem 1.7.

\demo{\bf Proof of Lemma 3.2}
Let $\beta_1\in C^\infty_0$ be supported in $(3/4,5/4)$ and equal to $1$
in $(7/8,9/8)$. Let $\chi$ be as in Lemma 3.1 supported 
in $\{|\xi_2|\in (4/5,6/5)\}$ and equal to $1$ in
$\{|\xi_2|\in (9/10,11/10)\}$.  Let
$$m_\k(\xi_1,\xi_3)\,=\,\beta(\xi_3)
\eta(4N(\xi_1/\xi_3-\alpha_\k))$$ and
$$\mu_{(N)}(\xi)\,=\,\sum_{\k=2}^N \chi(2^{-6\k}\xi_2)
m_\k(\xi_1,\xi_3).$$
In view of the properties of $\eta$, $\chi$, $\widetilde \chi$ and
the Marcinkiewicz multiplier theorem 
$$\|\mu_{(N)}\|_{M_{p2}}\,\le\,C_p
\|m_{(N)}\|_{M_{p2}},\qquad 1<p<\infty.
$$
Now assume $4\le p<\infty$. Let 
$$
R_\k\,=\,\{(x_1,x_3);\,|x_1-\alpha_\k x_3|\le 10^{-3}N,
|\alpha_\k x_1-x_3|\le 10^{-3}\}.
$$
For $\xi\in\supp m_\k$, $x\in R_\k$ we have $|x_1\xi_1+x_3\xi_3|\le \pi/4$ and therefore
$$
\Big|\int m_\k(\xi_1,\xi_3) e^{i(x_1\xi_1+x_3\xi_3)} d\xi_1d\xi_3\Big|
\ge
\Big|\int m_\k(\xi_1,\xi_3) \cos(x_1\xi_1+x_3\xi_3) d\xi_1d\xi_3\Big|
\ge c N^{-1}
$$
for some fixed positive constant $c$.
Let
$$
\align
\widetilde R_\k\,&=\,\{(x_1,x_3);\,10^{-4}N/2
\le|x_1-\alpha_\k x_3|\le 10^{-4}N,\,\,
10^{-4}/2\le|\alpha_\k x_1-x_3|\le 10^{-4}\}
\\
R_\k^*\,&=\,\{(x_1,x_3);\,
|x_1-\alpha_\k x_3|\le 10^{-4}N,\,\,
|\alpha_\k x_1-x_3|\le 10^{-4}\}
\endalign
$$
and let 
$\chi_\k$ be the characteristic function of $\widetilde R_\k$. Then
$$
\Cal F^{-1}[m_\k\Cal F\chi_\k]\ge c',\qquad x\in R_\k^*.
$$
By Lemma 3.1
$$
\Bigl\|\Bigl(\sum_{\k=2}^N\big|\Cal F_{\Bbb R^2}^{-1}[m_\k\Cal F_{\Bbb R^2}
 \chi_\k]\big|^2\Bigr)^{1/2}\Bigr
\|_{L^p(\Bbb R^2)}
\,\le\,C_{p} \|\mu_{(N)}\|_{M_{p2}}
\Bigl\|\Bigl(\sum_{\k=2}^N|\chi_\k|^2\Bigr)^{1/2}\Bigr\|_{L^p(\Bbb R^2)}
$$
Now
one verifies that
$$
\Bigl\|\Bigl(\sum_\k|\chi_\k|^2\Bigr)^{1/2}\Bigr\|_p\,\approx\,
N^{2/p}
$$
In view of the overlap of the rectangles $R_\k^* $ we
have for some small constant $c_1>0$, and for $|x|\le c_1 N$
and  for $|x|\le cN$ we have
$$
\Bigl(\sum_\k\big|\Cal F^{-1}[m_\k\Cal F \chi_\k](x)\big|^2\Bigr)^{1/2}\,
\approx\,N^{1/2}(1+|x|)^{-1/2}
$$
and consequently
$$\Bigl\|\Bigl(\sum_\k\bigl|\Cal F^{-1}[m_\k
\Cal F\chi_\k]\bigr|^2\Bigr)^{1/2}\Bigr\|_p\,\approx\,
\cases
N^{1/2}\qquad&\text{if } p>4
\\
N^{1/2}(\log N)^{1/4}\qquad&\text{if } p=4
\endcases
$$
This implies the assertion.\qed
\enddemo

Next we consider
 the class of homogeneous functions $m$ in 
$\Bbb R^3$ with the property
that the restrictions $h_\pm$ to the hyperplanes $\{\xi;\,\xi_3=\pm 1\}$
satisfy the hypotheses 
of the Marcinkiewicz multiplier theorem in the plane; that is
$$
\gathered
\|h\|_\infty\,\le\, A
\\
\sup_{j_1\in \Bbb N}\sup_{s_2} \int_{I_{j_1}}
\Big | \frac{\partial h}{\partial s_1}(s_1,s_2)\Big|\,ds_1\,\le\,A
\\
\sup_{j_2\in \Bbb N}\sup_{s_1} \int_{I_{j_2}}
\Big | \frac{\partial h}{\partial s_2}(s_1,s_2)\Big|\,ds_2\,\le\,A
\\
\sup_{j\in \Bbb N^2} \iint_{I_{j_1}\times I_{j_2}}
\Big | \frac{\partial^2 h}{\partial s_1\partial s_2}(s_1,s_2)\Big|\,ds_1\,ds_2
\le\,A
\endgathered
\tag 3.8
$$
where $I_{j_1}$ etc. is as in 
(1.1).
We  show that (3.8) is not sufficient to guarantee $m\in M_p$, for any
$p\neq 2$. The argument here follows Fefferman's solution \cite{16} of the 
multiplier problem for the ball (see also \cite{14}, \cite{21}).

Let $\frak a=\{\alpha_\k\}$ an arbitrary sequence of numbers in $[1,2)$ and let
$m_\k$ be defined in the first quadrant such that
$$
m_\k(\xi_1,\xi_3)\,=\,\cases
1,\qquad &1\le \xi_1/\xi_3\le\alpha_\k,\,
\\
0\qquad &\alpha_\k<\xi_1/\xi_3\le 2.
\endcases
\tag 3.9
$$
Let 
$$m_{\frak a}\,=\,\sum_\k\beta(2^\k\xi_2/\xi_3)m_\k(\xi_1,\xi_3).
$$
Suppose the assumptions (3.8) imply $m\in M_p$ for some $p\neq 2$. 
Then a limiting argument as in \cite{25, p.109} would imply that
$m_{\frak a}$ is an $L^p$ multiplier with norm independent of the choice of 
$\{\alpha_\k\}_{k\in\Bbb Z}$ and by (3.1) a corresponding statement on
$L^p(L^2)$ would follow.   However we have

\proclaim{\bf Lemma 3.3} 
The inequality
$$
\|\Cal F^{-1}[m_{\frak a}\Cal F f]\|_{L^p(L^2)}\,\le\, C
\|f\|_{L^p(L^2)}
$$
does not hold independently of $\frak a$ if $p\neq 2$.  
\endproclaim

For example   if we take
for $\frak a$ an
enumeration of the rationals in $[1,2)$ then $m_{\frak a}\in M_{p2}$ if 
and only if $p=2$.

\demo{\bf Proof}
Arguing as above the assumption
$m_{\frak a}\in M_p$ implies a vector-valued estimate for directional Hilbert 
transforms, namely
$$
\Bigl\|\Bigl(\sum_\k|H_\k f_\k|^2\Bigr)^{1/2}\Bigr\|_{L^p(\Bbb R^2)}
\,\le\, C
\Bigl\|\Bigl(\sum_\k| f_\k|^2\Bigr)^{1/2}\Bigr\|_{L^p(\Bbb R^2)}.
$$
where $H_\k$ is the Hilbert transform in the direction $(1,-\alpha_\k)$.
 But as in \cite{16} 
the existence of the Kakeya set prohibits such 
inequalities for  $p\neq 2$
 (unless further restrictions on the 
family of directions $(1,-\alpha_\k)$ are made). \qed
\enddemo

\bigskip

\head
{\bf 4.
Weighted norm inequalities in higher dimensions}
\endhead
We deduce Theorem 1.3 from  a weighted norm inequality; the 
procedure is analogous to
 Stein's proof of the 
H\"ormander multiplier theorem (see \cite{25, ch.IV}). Here
however the  positive operator which controls the problem is 
not the Hardy-Littlewood maximal 
operator but a multiple  iteration of 
variants of Kakeya-type maximal operators. The main step  of the argument
is contained in 
Lemma 4.7; one  proves a weighted
inequality for a variant of
C\'ordoba's sectorial square-function.

For $i=1,\dots,d-1$ let 
 $\Cal R^{(i,d)}_{n_i}$ be the family of all rectangles with dimensions 
$1\times 2^{n_i}$, centered at the origin in the $x_i-x_d$ plane 
and let
$\widetilde{\Cal R}^{(i,d)}_{n_i}$ be the family of all parallelograms
of the form
$\{(x_i,x_d); (2^{k_i}x_i, 2^{k_d}x_d)\in R_0\}$  where $R_0\in
\Cal R^{(i,d)}_{n_i}$ and $k_i,\,k_d$ are integers. Let
$$
M^{(i,d)}_{n_i}\om(x_1,\dots,x_d)\,=\,
\sup_{R\in\widetilde{\Cal R}^{(i,d)}_{n_i}}\frac 1{|R|}\int_R
|\om(x_1,\dots,x_{i-1},x_i-y_i,x_{i+1},\dots,x_{d-1},x_d-y_d)|\,dy_idy_d.
$$
$M_{n_i}^{(i,d)}$ is a  variant of the  Kakeya maximal
function, invariant under the dilations
$(x_i,x_d)\mapsto (2^{k_i}x_i,2^{k_d}x_d)$.
The proof of the $L^2$-estimate 
in \cite {10} can be  easily modified to yield
$$\|M_{n_i}^{(i,d)}\om\|_2\,\le\, C n_i\|\om\|_2;$$
for a more singular variant see also \cite{11}.

Next, for $n=(n_1,\dots,n_{d-1})$
define 
$$\Cal M_n\,=\, M^{(d-1,d)}_{n_{d-1}}
\circ M^{(d-2,d)}_{n_{d-2}}\circ\dots\circ M^{(1,d)}_{n_1}$$
and, for $N\in\Bbb N$, let $\Cal M^N_n=\Cal M_n\circ\dots\circ\Cal M_n$ be the
$N$-fold application of the operator $\Cal M_n$.  Finally, if 
$M_{(i)}$ denotes the Hardy-Littlewood maximal operator with
respect to the variable $x_i$ let
$$
\widetilde {\Cal M}^N_n\,=\,
M_{(1)}\circ\dots\circ M_{(d)}\circ
\Cal M^N_n\circ
M_{(1)}\circ\dots\circ M_{(d)}.
$$

\proclaim{Theorem 4.1} Let $\gamma>1/2$ and suppose that
$$\sup_{t\in(\Bbb R_{+})^{d-1}}
\|\beta_{(1)}\otimes\dots\otimes\beta_{(d-1)}
\, g_\pm(t_1\cdot,\dots,t_{d-1}\cdot)
\|_{\Cal H^2_\gamma(\Bbb R^{d-1})}\,\le\, B_\gamma\,<\,\infty
\tag  4.1
$$
Let $m$ be the homogeneous extension of $g_\pm$ and define $T$
by $\widehat{ Tf}(\xi)=m(\xi)\widehat f(\xi)$.
Let $0<\ep<\gamma-1/2$, let $N(\ep)$ be the smallest positive integer 
$> 3+2/\ep$ and define $\frak M_\ep$ by
$$
\frak M_\ep \om\,=\,\sum_{n\in\Bbb N_0^{d-1}}
2^{-\ep(n_1+\dots+n_{d-1})}
\widetilde{\Cal M}_n^{N(\ep)}\om.
$$
Then for $s>1$
$$
\int|Tf(x)|^2\om(x)\,dx\,\le\,
C_{\ep, s}B_\gamma
\int|f(x)|^2 (\frak M_\ep[\om^s])^{1/s}dx.
\tag 4.2
$$
\endproclaim

\demo{\bf Proof of Theorem 1.5}
Since the  operator 
$\om\mapsto (\frak M_\ep(|\om|^s))^{1/s}$ is bounded on $L^q$, $q>2s/(1-\ep)$,
the weighted norm inequality (4.2) and  duality imply  under the 
assumption (4.1) 
that $T$ is bounded on $L^p$, for $2\le p\le 4. $
The general result of Theorem 1.5 follows then by interpolation, using the
technique in \cite {9}.\qed

\enddemo

Before we prove Theorem 4.1 we recall a few facts about 
vector-valued weighted norm inequalities. First if  $H_1$, $H_2$ are 
Hilbert spaces and
 $\Cal K$ is a
convolution kernel in $\Bbb R$,
 with values in the space $\Cal B(H_1,H_2)$  of bounded 
operators, then $\Cal K$ is called a regular singular integral operator if
$$
\align
|\widehat{ \Cal K}(\xi)|_{\Cal B(H_1,H_2)}&\le C
\\
|\Cal K(x)|_{\Cal B(H_1,H_2)}&\le C|x|^{-1}
\\
|\Cal K(x-y)-\Cal K(x)|_{\Cal B(H_1,H_2)}&\le C|y|^\delta|x|^{-1-\delta},\qquad
|x|>2|y|>0;
\endalign
$$
here $0<\delta\le 1$ is fixed.
By a  vector-valued version of a theorem of C\'ordoba and Fefferman 
(\cite{13}, see also \cite{18, ch.IV.3})
there is an inequality
$$
\int|\Cal K*f(x)|^p_{H_2}\om(x) dx\,\le\,
C_{\sigma, p}\int|f(x)|_{H_1}^p(M(|\om|^\sigma))^{1/\sigma}(x)\, dx
\tag 4.3
$$
where $1<p<\infty$,
$\sigma>1$.

Littlewood-Paley functions can be associated with regular singular 
integral operators.
Let $\beta\in C^\infty_0(1/2,2)$ then it is straightforward to check that
the operator $\{f_\k\}_{\k\in \Bbb Z}\mapsto \sum \Cal F^{-1}[\beta(2^\k\cdot)
\Cal F f]$ is a $\Cal B(\ell^2,\Bbb R)$-valued 
regular singular integral operator. Likewise 
the adjoint operator
$f\mapsto \{\Cal F^{-1}[\beta(2^\k\cdot)\Cal F f]\}_{\k\in\Bbb Z}$ is a 
$\Cal B(\Bbb R,\ell^2)$-valued  regular singular integral operator.
Here $\ell^2$ may refer to a space of sequences with values in a Hilbert space.

Next let  $k\in \Bbb Z^d$  and denote by
 $\L_k$ be the standard Littlewood-Paley operator
with multiplier $\prod_{i=1}^d\beta(2^{k_i}\xi_i)$.
Then a repeated application of (4.3) yields

\proclaim{Lemma 4.2}
For $1<p<\infty$,
$s>1$ we have the inequalities
$$
\gather
\int\Big|\sum_{k\in\Bbb Z^d} \L_k f_k\Big|^p\om(x) dx\,\le\,
C_{s, p}\int\Bigl(\sum_{k\in\Bbb Z^d}|f_k(x)|^2\Bigr)^{p/2}
(M_{(1)}\circ\dots\circ M_{(d)}[\om^s])^{1/s}(x)\, dx
\\
\int\Bigl(\sum_{k\in\Bbb Z^d}|\L_k f|^2\Bigr)^{p/2}\om(x) dx\,\le\,
C_{s, p}\int|f(x)|^p
(M_{(1)}\circ\dots\circ M_{(d)}[\om^s])^{1/s}(x)\, dx.
\endgather
$$
\endproclaim

We need also the following pointwise estimate  concerning a 
square-function involving translates of a fixed 
Schwartz-function $\eta$. It implies  $L^p$-boundedness  for
$p>2$, a result which is due to Carleson. 
A proof of the  pointwise estimate can be found 
in \cite{24}, see also \cite{12}, \cite{18}.

\proclaim{Lemma 4.3} Let $\eta$ be a Schwartz function in $\Bbb R^d$ 
and let $A\in GL(d,\Bbb R)$.
Then
$$
\sum_{k\in \Bbb Z^d}\bigl|\Cal F^{-1}[\eta(A\cdot-k)\Cal F f](x)]
\bigr|^2
\,\le\,
C_N\int\frac{|f(x-A^t y)|^2}{(1+|y|)^N}\,dy.
\tag 4.4
$$
\endproclaim

\subheading{\bf Proof of Theorem 4.1}

There is no loss of generality in assuming that $m$ is supported in
$\{\xi;\,\xi_i\ge 0,\,i=1,\dots,d\}$. Setting
$$\phi(\xi)=\prod_{i=1}^d\beta(\xi_i)
\tag 4.5
$$
 we decompose
$$
m(\xi)\,=\,\sum_{k\in \Bbb Z^d}m_k(2^{k_1}\xi_1,\dots,2^{k_d}\xi_d)
$$
where
$$
m_k(\xi)\,=\,\phi(\xi) g_k(\xi_1/\xi_d,\dots,\xi_{d-1}/\xi_d)
$$
and $g_k$ has compact support in $(1/2,2)^{d-1}$.
Note that $g_k=g_{k'}$ if $k_i-k_d=k_i'-k_d'$,
$i=1,\dots,d-1$. We introduce a further decomposition
using the dyadic smooth cutoff functions 
$\psi_n=\psi_{n_1}\otimes\dots\otimes\psi_{n_{d-1}}$ 
({\it cf.} the second definition of the space $\Cal H^q_\alpha$ 
in the introduction).
We decompose
$$
\align
m_k(\xi)\,&=\,\sum_{n\in(\Bbb N_0)^{d-1}}
\phi(\xi)
\,g_k*\widehat{\psi_n}(\xi_1/\xi_d,\dots,\xi_{d-1}/\xi_d)
\\
\,&=\,\sum_{n\in(\Bbb N_0)^{d-1}}
\phi(\xi)m_k^n(\xi).
\tag 4.6
\endalign
$$
We may write
$$g_k*\widehat{\psi_n}\,=\,
g_k^n*\widehat{\psi_n}
\tag 4.7
$$
where
$$g^n_k\,=\,g_k*\widehat{\widetilde\psi_n}$$
and $\widetilde \psi_n=\widetilde\psi_{n_1}\otimes\dots\otimes
\widetilde\psi_{n_{d-1}}$ 
is similarly defined as $\psi_n$ (say, with 
$\widetilde\psi_{n_i}$ supported in 
$\pm[2^{n_i-2},2^{n_i+2}]$,
 equal to $1$
in $\supp\,\psi_{n_i}$). 
Let us note in passing that in view of the support properties of the 
Fourier transform  of $g^n_k$ we have the following version of 
Sobolev's  imbedding theorem
$$
\sup_{s_{d_1+1},\dots,s_{d-1}}\|g^n_k(\cdot,s_{d_1+1},\dots,s_{d-1})
\|_{L^p(\Bbb R^{d_1})}\,\le\,
C2^{(n_{d_1+1}+\dots+n_{d-1})/p}\|g^n_k\|_{L^p(\Bbb R^{d-1})},
\tag 4.8
$$
see the argument in \cite{27, p.18}.

Let $T_k^n$ be defined by
$$
\widehat {T^n_k f}(\xi)\,=\,[\phi m^n_k](2^{k_1}\xi_1,\dots,2^{k_d}\xi_d)
\widehat f(\xi).
\tag 4.9
$$

Let $0<\ep'<\ep$, say $\ep'=\ep/2$.  An application of Lemma 4.2 
shows that it suffices to prove the 
inequality
$$\multline
\sum_{k\in\Bbb Z^d}\int|T^n_k \L_kf(x)|^2\om(x)\,dx\,
\\
\le\,
C_{N} 2^{(n_1+\dots+n_{d-1})(\frac 12+\ep')}\|g^n_k\|_2\sum_{k\in\Bbb Z^d}
\int|\L_kf(x)|^2 \Cal M^N_n\om(x)\, dx,\qquad N>2+\frac 1{\ep'}.
\endmultline
\tag 4.10
$$
In order to avoid complicated notation we shall assume $d=3$ in what follows.
This case is entirely 
typical of the general situation in higher dimensions.

In order to use the homogeneity of 
the multipliers we have to introduce finer decompositions of $g^k_n$.
For $\nu_1= 2^{n_1-3}, 2^{n_1-3}+1, \dots, 2^{n_1+3}$ and 
$\nu_2= 2^{n_2-3}, 2^{n_2-3}+1, \dots, 2^{n_2+3}$ let 
$$u_\nu\,=\,
(u_{\nu_1}^1,u_{\nu_2}^2)\,=\,(2^{-n_1}\nu_1,2^{-n_2}\nu_2)
\tag 4.11
$$
 and 
$$
I_\nu\,=\,I_{\nu_1}^1\times I_{\nu_2}^2\,=\,
[u_{\nu_1}^1-2^{-n_1-1},u_{\nu_1}^1+2^{-n_1-1}]\times
[u_{\nu_2}^2-2^{-n_2-1},u_{\nu_2}^2+2^{-n_2-1}].
$$
Furthermore let
$$
\align
{}^cI^1\,&=\,\Bbb R\setminus\cup_{\nu_1=2^{n_1-3}}^{2^{n_1+3}}I^1_{\nu_1}
\,=\,\Bbb R\setminus[\frac 18 -2^{-n_1-1},\,8+2^{-n_1-1}]
\\
{}^cI^2\,&=\,\Bbb R\setminus\cup_{\nu_2=2^{n_2-3}}^{2^{n_2+3}}I^2_{\nu_2}
\,=\,\Bbb R\setminus[\frac 18 -2^{-n_2-1},\,8+2^{-n_2-1}].
\endalign
$$
Setting
$$
\aligned
g_{k\nu}^n(s)\,&=\,\int_{I_\nu} g_k^n(u)
\widehat{\psi_n} (s-u)\, du 
\\
g_{k\nu_1}^{n,1}(s)\,&=\,
\int_{ I^1_{\nu_1}\times {}^cI^2}
g_k^n(u)\widehat{\psi_n} (s-u)\, du 
\\
g_{k\nu_2}^{n,2}(s)\,&=\,
\int_{{}^cI^1\times I^{2}_{\nu_2}}
g_k^n(u)\widehat{\psi_n} (s-u)\, du 
\\
\rho_k^n(s)\,&=\,
\int_{ {}^cI^1\times{}^cI^2}
g_k^n(u)\widehat{\widetilde \psi_n} (s-u)\, du.
\endaligned
\tag 4.12$$
we split
$$
T_k^n\,=\,
\sum_{\nu_1=2^{n_1-3}}^{2^{n_1+3}}
\sum_{\nu_2=2^{n_2-3}}^{2^{n_2+3}}
T_{k\nu}^n
+\sum_{\nu_1=2^{n_1-3}}^{2^{n_1+3}}
T_{k\nu_1}^{n,1}
+\sum_{\nu_2=2^{n_2-3}}^{2^{n_2+3}}
T_{k\nu_2}^{n,2}
+T_k^{n,0}
\tag 4.13
$$
where
$$
\aligned
\widehat{T^n_{k\nu}f}(\xi)\,&=\,
g_{k\nu}^{n}(2^{k_1-k_3}\xi_1/\xi_3,2^{k_2-k_3}\xi_2/\xi_3)
\phi(2^{k_1}\xi_1,2^{k_2}\xi_2,2^{k_3}\xi_3)
\,\widehat f(\xi)
\\
\widehat{T^{n,1}_{k\nu_1}f}(\xi)\,&=\,
g_{k\nu_1}^{n,1}(2^{k_1-k_3}\xi_1/\xi_3,2^{k_2-k_3}\xi_2/\xi_3)
\phi(2^{k_1}\xi_1,2^{k_2}\xi_2,2^{k_3}\xi_3)\,\widehat f(\xi)
\\
\widehat{T^{n,2}_{k\nu_2}f}(\xi)\,&=\,
g_{k\nu_2}^{n,2}(2^{k_1-k_3}\xi_1/\xi_3,2^{k_2-k_3}\xi_2/\xi_3)
\phi(2^{k_1}\xi_1,2^{k_2}\xi_2,2^{k_3}\xi_3)\widehat f(\xi)
\\
\widehat{T^{n,0}_kf}(\xi)\,&=\,
\rho_k^n(2^{k_1-k_3}\xi_1/\xi_3,2^{k_2-k_3}\xi_2/\xi_3)
\phi(2^{k_1}\xi_1,2^{k_2}\xi_2,2^{k_3}\xi_3)\,\widehat f(\xi).
\endaligned
\tag 4.14$$
We set
$$
\aligned
b^n_{k\nu}\,&=\,\sup_{u\in I_\nu} |g^k_n(u)|
\\
b^{n,1}_{k\nu_1}\,&=\,\sup_{u\in I^1_{\nu_1}\times {}^cI^2}
|g^k_n(u)|
\\
b^{n,2}_{k\nu_2}\,&=\,
\sup_{u\in{}^cI^1\times I^{2}_{\nu_2}}
|g^k_n(u)|
\\
b^n_{k}\,&=\,\sup_{u\in {}^cI^1\times{}^cI^2} |g^k_n(u)|.
\endaligned
\tag 4.15
$$
Since the Fourier transform of $g^n_k$ is supported in 
$[-2^{n_1+3},2^{n_1+3}]\times
[-2^{n_2+3},2^{n_2+3}]$, 
suitable  variants of the Plancherel-Polya theorem (see \cite{27, p.19})
and the Sobolev embedding theorem imply
$$
\Bigl(
\sum_{\nu_1=2^{n_1-3}}^{2^{n_1+3}}
\sum_{\nu_2=2^{n_2-3}}^{2^{n_2+3}}
[b_{k\nu}^n]^r\Bigr)^{1/r}
\,\le\,C_r\,2^{(n_1+n_2)/r}\|g^n_k\|_r,\qquad 0<r\le\infty
\tag 4.16
$$
with the appropriate interpretation for $r=\infty$; moreover we have
$$
\align
\Bigl(\sum_{\nu_1=2^{n_1-3}}^{2^{n_1+3}}
[b_{k\nu_1}^{n,1}]^r\Bigr)^{1/r}
\,&\le\,C\,2^{n_1/r}\sup_{s_2}\|g^n_k(\cdot,s_2)\|_{L^r(\Bbb R)}
\\
&\le\,C\,2^{(n_1+n_2)/r}\|g^n_k\|_{L^r(\Bbb R^2)}
\tag 4.17
\endalign
$$
and a similar statement with the $s_1$ and $s_2$ variables interchanged.
Also by (4.8) $b^n_k$ is bounded by $C2^{(n_1+n_2)/r}\|g^n_k\|_r$.

We  need pointwise estimates for the  convolution kernels
$K^n_{k\nu}$, $K^{n,1}_{k\nu_1}$, $K^{n,2}_{k\nu_2}$, $K^{n,0}_k$ 
of the operators
$T^n_{k\nu}$, $T^{n,1}_{k\nu_1}$, $T^{n,2}_{k\nu_2}$, $T^{n,0}_k$, 
respectively.

\proclaim{\bf Lemma 4.4} 
Let  $e_\nu=(u_{\nu_1}^1,u_{\nu_2}^2,1)$,
$e^1_{\nu,1}=(u_{\nu_1}^1,0,1)$ and
$e^2_{\nu,2}=(0,u_{\nu_2}^2,1)$.
Let
$$
\align
W^n_{\nu N}(x)\,&=\,2^{-n_1-n_2}
(1+|\inn{e_\nu}{x}|)^{-N}(1+2^{-n_1}|x_1|)^{-N}
(1+2^{-n_2}|x_2|)^{-N}
\\
W^{n,1}_{\nu_1 N}(x)\,&=\,2^{-n_1}
(1+|\inn{e^1_{\nu_1}}{x}|)^{-N}(1+2^{-n_1}|x_1|)^{-N}
(1+|x_2|)^{-N}
\\
W^{n,2}_{\nu_2 N}(x)\,&=\,2^{-n_2}
(1+|\inn{e^2_{\nu_2}}{x}|)^{-N}
(1+|x_1|)^{-N}(1+2^{-n_2}|x_2|)^{-N}
\\
W^{n,0}_{N}(x)\,&=\,(1+|x_1|)^{-N}(1+|x_2|)^{-N}
(1+|x_3|)^{-N}.
\endalign
$$
Let
$$U^n_{k\nu,N}(x)\,=\,2^{-k_1-k_2-k_3}W^n_{k\nu}
(2^{-k_1}x_1,2^{-k_2}x_2,2^{-k_3}x_3)
$$
and similarly define  $U^{n,1}_{k\nu_1,N}$,
$U^{n,2}_{k\nu_2,N}$, $U^{n,0}_{k,N}$.
Then
$$
\aligned
|\partial_x^\gamma K^n_{k\nu}(x)|
\,&\le\,
C_{\gamma N} b^n_{k\nu}
 2^{-k_1\gamma_1-k_2\gamma_2-k_3\gamma_3} U^n_{k\nu, N}(x)
\\
|\partial_x^\gamma K^{n,1}_{k\nu_1}(x)|
\,&\le\,
C_{\gamma N} b^{n,1}_{k\nu_1} 2^{-n_2N}
 2^{-k_1\gamma_1-k_2\gamma_2-k_3\gamma_3} U^{n,1}_{k\nu_1, N}(x)
\\
|\partial_x^\gamma K^{n,2}_{k\nu_2}(x)|
\,&\le\,
C_{\gamma N} b^{n,2}_{k\nu_2} 2^{-n_1N}
2^{-k_1\gamma_1-k_2\gamma_2-k_3\gamma_3} U^{n,2}_{k\nu_2, N}(x)
\\
|\partial_x^\gamma K^{n,0}_{k}(x)|
\,&\le\,
C_{\gamma N} b^n_k
2^{-n_1N}2^{-n_2N}
2^{-k_1\gamma_1-k_2\gamma_2-k_3\gamma_3} U^{n,0}_{k, N}(x).
\endaligned
\tag 4.18
$$
\endproclaim
\demo{\bf Proof}
First consider $K^n_{k\nu}$.
Using the homogeneity of the  multiplier and the decay properties of
$\widehat{\psi_n}$ we see that
$$
\multline
\Bigl|\partial_{\xi_1}^{N_1}\partial_{\xi_2}^{N_2}
{\inn{e_\nu}{\nabla_\xi}}^{N_3}
\Bigl[
\phi(\xi)
\int_{I_\nu}g^n_{k}(u)
\widehat{\psi_n}(\frac{\xi_1}{\xi_3}-u_1,\frac{\xi_2}{\xi_3}-u_2)du 
\Bigr]\Bigr|\,
\\
\,\le\,
C(N_1,N_2,N_3,M) \,b^n_{k\nu}\, 
\frac{2^{N_1n_1}}{(1+2^{n_1}|\xi_1/\xi_3-u_{\nu_1}^1|)^{M}}
\frac{2^{N_2 n_2}}{(1+2^{n_2}|\xi_2/\xi_3-u_{\nu_2}^2|)^{M}}.
\endmultline
\tag 4.19
$$
Using integration by parts we obtain 
$$
\align
2^{k_1+k_2+k_3}&|K^n_{k\nu}(2^{k_1}x_1,2^{k_2}x_2,2^{k_3}x_3)|\\
\,\le\,&C_{NM}\int_{[1/2,2]^3}
(1+2^{n_1}|\xi_1/\xi_3-u_{\nu_1}^1|)^{-M}
(1+2^{n_2}|\xi_2/\xi_3-u_{\nu_2}^2|)^{-M} d\xi
\, 2^{n_1+n_2}W^n_{\nu N}(x)
\\
\le\, &C_{N} \, W^n_{\nu N}(x).
\endalign
$$
In view of the compact support of $\phi$ we get the same estimates
for the derivatives of
the left hand side 
and the desired estimates for $K^n_{k\nu}$ and its derivatives follow.

The estimate for $K^{n,0}_k$ has nothing to do with homogeneity: By the decay 
properties of $\widehat{\psi_n}$ 
we have
$$
\big|\partial_{s_1}^{\gamma_1}\partial_{s_2}^{\gamma_2}\rho_k^n(s)\big|\,\le\,
C_{\gamma N}2^{-(n_1+n_2)N}(1+|s_1|)^{-N}(1+|s_2|)^{-N}
$$
and hence
$$
\big|\partial_{\xi_1}^{\gamma_1}\partial_{\xi_2}^{\gamma_2}
\partial_{\xi_3}^{\gamma_3}
[\phi(\xi)\rho_k^n(\xi_1/\xi_3,\xi_2/\xi_3)]\big|
\,\le\,C_{\gamma N}2^{-(n_1+n_2)N}.
$$
The desired estimate for $K^{n,0}_{k}$ follows by integration by parts.
In the proof of the estimate for $K^{n,1}_k$ 
we replace (4.19) by
$$
\big|\partial_{\xi_1}^{N_1}\partial_{\xi_2}^{N_2}
{\inn{e_\nu^1}{\nabla_\xi}}^{N_3}
[\phi(\xi)
g_{k\nu_1}^{n,1}(\xi_1/\xi_3,\xi_2/\xi_3)]\big|\,
\,\le\,
C(N_1,N_2,N_3,M) \,b^n_{k\nu}\, 
\frac{2^{N_1n_1}2^{(N_2+N_3-M)n_2}}
{(1+2^{n_1}|\xi_1/\xi_3-u_{\nu_1}^1|)^{M}}
$$
and argue as above.\qed
\enddemo

In what follows we  shall  denote by 
$\tbeta$ a function which is 
similar to $\beta$ but equals $1$ on the support of $\beta$.
Next let $\chi\in C^\infty_0(\Bbb R)$ be supported in $(3/4,3/4)$ such that
$\sum_{\k\in \Bbb Z} \chi^2(\cdot-\k)\equiv 1$.
Again let $\widetilde \chi\in C^\infty_0$ be defined similarly to $ \chi$ but
equal to $1$ on the support of $\chi$.
We define the operator $A^n_{k\nu}$ by
$$
\widehat{ A^n_{k\nu}f}(\xi)\,=\,
\chi(2^{n_1}(2^{k_1}\xi_1-u^1_{\nu_1} 2^{k_3}\xi_3))
\chi(2^{n_2}(2^{k_2}\xi_2-u^2_{\nu_2} 2^{k_3}\xi_3))
\tbeta^2(2^{k_3}\xi_3)\widehat f(\xi).
$$

\proclaim{\bf Lemma 4.5}
There is a weighted norm inequality
$$\int\sum_k|T^n_k \L_kf(x)|^2\om(x)\,dx\,\le\,
C 2^{n_1+n_2}\|g^n_k\|_2^2\,\int
\sum_{k,\nu}|A^n_{k\nu}\L_kf(y)|^2 \Cal M_n\om(y)\, dy.
\tag 4.20
$$
\endproclaim

\demo{\bf Proof}
Set 
$ S^n_{k\nu\mu}
=T^n_{k\mu} A^n_{k\nu}$,
$ S^{n,1}_{k\nu\mu_1}
=T^{n,1}_{k\mu_1} A^n_{k\nu}$,
$ S^{n,2}_{k\nu\mu_2}
=T^{n,2}_{k\mu_2} A^n_{k\nu}$ and
$ S^{n,0}_{k\nu}=T^{n,0}_k A^n_{k\nu}$.
Then
$$
T^n_k \L_kf\,=\,
\sum_\nu\bigl[
\sum_\mu S^n_{k\nu\mu}+\sum_{\mu_1} S^{n,1}_{k\nu\mu_1}
+\sum_{\mu_2} S^{n,2}_{k\nu\mu_2}+S^{n,0}_{k\nu}\bigr]
A^n_{k\nu} \L_kf.
\tag 4.21
$$
Let $H^n_{k\nu\mu}$, 
$H^{n,1}_{k\nu\mu_1}$, 
$H^{n,2}_{k\nu\mu_2}$ and
$H^{n,0}_{k\nu}$
 be the convolution kernels of the operators
$S^n_{k\nu\mu}$, 
$S^{n,1}_{k\nu\mu_1}$, 
$S^{n,2}_{k\nu\mu_2}$ and
$S^{n,0}_{k\nu}$, respectively.
Fix $N$ (say equal to $100$) and let $U^n_{k\nu}\equiv U^n_{k\nu, 100}$ etc.
The proof of Lemma 4.4 shows that
$$
\align
|H^n_{k\nu\mu}(x)|\,&\le\,
C\frac{b^n_{k\mu}}{(1+|\mu_1-\nu_1|^2)(1+|\mu_2-\nu_2|^2)}
U^n_{k\nu}(x)
\\
|H^{n,1}_{k\nu\mu_1}(x)|\,&\le\,
C_N 2^{-n_2N} \frac{b^{n,1}_{k\mu_1}}{(1+|\mu_1-\nu_1|^2)}
U^{n}_{k\nu}(x)
\\
|H^{n,2}_{k\nu\mu_2}(x)|\,&\le\,
C_N 2^{-n_2N} \frac{b^{n,2}_{k\mu_2}}{(1+|\mu_2-\nu_2|^2)}
U^{n}_{k\nu}(x)
\\
|H^{n,0}_{k\nu}(x)|\,&\le\,
C_N 2^{-n_1 N}2^{-n_2N} b^{n}_{k}
U^{n}_{k\nu}(x).
\endalign
$$
We observe that
$\|U^n_{k\nu}\|_1\,\le\, C$ is bounded
uniformly in $\nu$, $n$, $k$.
Therefore
$$
\align
&\Big|\sum_\nu \sum_\mu S^n_{k\nu\mu} A^n_{k\nu}\L_k f(x)\Big|
\\
\,\le\,&C\sum_\nu
\sum_\mu\frac{|b^n_{k\mu}|}{(1+|\mu_1-\nu_1|^2)
(1+|\mu_2-\nu_2|)^2}
\Bigl(\int U^n_{k\nu}(x-z)dz\Bigr)^{1/2}
\Bigl(\int|A^n_{k\nu}\L_kf(y)|^2
U^n_{k\nu}(x-y)dy\Bigr)^{1/2}
\\
\le \,&C
\Bigl(\sum_\nu\Big|
\sum_\mu\frac{|b^n_{k\mu}|}{(1+|\mu_1-\nu_1|^2)
(1+|\mu_2-\nu_2|)^2}
\Big|^2\Bigr)^{1/2}
\Bigl(\sum_\nu\int|A^n_{k\nu}\L_kf(y)|^2
U^n_{k\nu}(x-y)dy\Bigr)^{1/2}
\\
\le\,&C
\Bigl(\sum_\mu|b^n_{k\mu}|^2\Bigr)^{1/2}
\Bigl(\sum_\nu\int|A^n_{k\nu}\L_kf(y)|^2
U^n_{k\nu}(x-y)dy\Bigr)^{1/2}
\\
\le\,&C 2^{n_1+n_2} \|g^n_k\|_2
\Bigl(\sum_\nu\int|A^n_{k\nu}\L_kf(y)|^2
U^n_{k\nu}(x-y)dy\Bigr)^{1/2}
\endalign
$$ 
where for the last inequality we have used (4.16).
Using also (4.17)  we derive the same inequality for the 
other three remaining terms in (4.21) and obtain
$$\Bigl(\sum_k|T^n_k \L_kf(x)|^2\Bigr)^{1/2}\,\le\,
C 2^{(n_1+n_2)/2}\|g^n_k\|_2
\Bigl(\int\sum_{k,\nu}|A^n_{k\nu}\L_kf(y)|^2
U^n_{k\nu}(x-y)\,dy
\Bigr)^{1/2}.
\tag 4.22
$$ 
Finally there is the pointwise estimate
$$
\sup_{k\nu}  U_{k\nu}^n*|\om|(x)\,\le\, C \Cal M_n\om(x)
\tag 4.23
$$
and (4.22) and (4.23) imply (4.20).\qed
\enddemo

\proclaim{\bf Proposition 4.6}
There is the weighted norm inequality
$$
\sum_{k,\nu}\int
|A^{n}_{k\nu}f_k(x)|^2 \om(x)\, dx\,\le\,
C_{N_0} 2^{(n_1+n_2)2\epsilon'}
 \sum_{k}\int |f_k(x)|^2\Cal M^{N_0}_n \om(x)\,dx,\qquad N_0>2+1/\ep'.
\tag 4.24
$$
\endproclaim

\demo{\bf Proof}
It is convenient to introduce a
 decomposition in the $\xi_3$-variable which will give the 
factor of $2^{2\ep'(n_1+n_2)}$.
We define for $(\la_1,\la_2)\in\Bbb Z^2$ operators 
$V^{n\ep'}_{k\la}$ by
$$\widehat{V^{n\ep'}_{k\la} f}(\xi)\,=\,
\chi^2(2^{k_3}\xi_3-2^{n_1\ep'}\la_1)
\chi^2(2^{k_3}\xi_3-2^{n_2\ep'}\la_2)
\widehat f(\xi)
$$
Observe that
 $A^{n}_{k\nu}$ is a sum of no more then $O(2^{\ep'(n_1+n_2)})$ operators
$V^{n\ep'}_\la A^{n}_{k\nu}$  where 
$2^{n_1\ep'}\la_1\in(1/20,20)$ and
$2^{n_2\ep'}\la_2\in(1/20,20)$.
Therefore it suffices to show that for those $\la$ the inequality
$$
\sum_{k,\nu}\int
|V^{n\ep'}_{k\la}
A^{n}_{k\nu}f_k(x)|^2 \om(x)\, dx\,\le\,
C_{N_0} \sum_{k}\int |f_k(x)|^2\Cal M^{N_0}_n \om(x)\,dx,\qquad N_0>2+1/\ep',
\tag  4.25
$$
holds.
In order to show (4.25) we first prove an inequality for an analogous problem
in two dimensions.

\proclaim{\bf Lemma 4.7}
Let $\delta\ll 1$ and let  $m,\mu,\rho$ be  integers such that
$m>0$,
$2^{-m}\mu\in(1/20,20)$ and
$2^{-m\delta}\rho\in(1/20,20)$.
Let $B^m_{\mu}$, $C^{m\delta}_\rho$
 be the operators acting on functions in $\Bbb R^2$ 
defined by
$$
\align
\widehat{ B^m_{\mu}f}(\xi)\,&=\,
\chi(2^{m}(\xi_1-2^{-m}\mu\xi_2))
\widehat f(\xi)
\\
\widehat{C^{m\delta}_\rho f}(\xi)\,&=\,
\chi^2(2^{m\delta}(\xi_2-2^{-m\delta}\rho))
\widehat f(\xi)
\endalign
$$
Let $l\le \max\{1, m\delta\}$.
Then
$$
\sum_{\mu}\int
|B^{m}_{\mu}C^{m\delta}_\rho f(x)|^2 \om(x)\, dx\,\le\,
C \sum_{\mu'}\int |B^{m-l}_{\mu'}
C^{m\delta}_\rho f(x)|^2 M^{(1,2)}_m \om(x)\,dx.
$$
\endproclaim

\demo{\bf Proof}
Let 
$$\align
R_{\mu\rho}\,&=\,\{\xi;\,|\xi_1-2^{-m}\mu\xi_2|\le 2^{-m+1};
|\xi_2-2^{-l}\rho|\le 2^{-l+1}\}
\\
\widetilde R_{\mu\rho}\,&=\,\{\xi;\,|\xi_1-2^{-m}\mu\xi_2|\le 2^{-m+2};
|\xi_2-2^{-l}\rho|\le 2^{-l+4}\}.
\endalign
$$
Let $\xi'\in R_{\mu'\rho}$ and suppose that $|\mu-\mu'|\le 2^{-l+2}$.
Let $a_{\mu-\mu'}=(2^{-m}(\mu-\mu'),0)$. Then
$\xi'-a_{\mu-\mu'}\in \widetilde R_{\mu\rho}$.
Thus 
$$
R_{\mu'\rho}\subset a_{\mu-\mu'}+\widetilde R_{\mu\rho}
$$
Define 
$$\widehat{\Gamma_{\mu\mu'\rho}^{ml} f}(\xi)\,=\,
\chi(2^{m-10}(\xi_1-2^{-m}(\mu-\mu')-2^{-m}\mu\xi_2))
\widetilde\chi(2^{m\delta}(\xi_2-2^{-m\delta}\rho))\widehat f(\xi).
$$
and define $\widetilde {C^{m\delta}_\rho}$ 
similarly as 
$C^{m\delta}_\rho$ 
(with $\chi$ replaced by $\widetilde \chi$).
Then
$$
\sum_{\mu}|B^{m}_{\mu}C^{m\delta}_\rho 
f(x)|^2 \,\le\,C\,
\sum_{\mu'}\sum_{\mu}
|\widetilde C^{m\delta}_\rho B^{m}_{\mu}
 \Gamma_{\mu\mu'\rho}^{ml}
C^{m\delta}_\rho B^{m-l}_{\mu'}f(x)|^2.
$$
An integration by parts argument shows that the convolution kernel of
$\widetilde C^{m\delta}_\rho B^{m}_{\mu}$ is bounded by $C_N$ times
$$
 w^{m,\delta}_{\mu N}(x)\,=\,\frac{2^m}{(1+2^m|\inn{x}{e_\mu}|)^{N}}
\frac{2^{m\delta}}{(1+2^{m\delta}|\inn{x}{e_\mu^\perp}|)^{N}}
$$
if $e_\mu=(1,-2^{-m}\mu)$, $e_\mu^\perp=(2^{-m}\mu,1)$ and if 
$2^{-m}\mu\approx 1$. Now the argument which lead to (4.22) and Lemma 4.3
show that for fixed $\mu'$
$$
\align
\sum_{|\mu-\mu'|\le l}
\int
&|\widetilde C^{m\delta}_\rho B^{m}_{\mu}
 \Gamma_{\mu\mu'\rho}^{ml}
C^{m\delta}_\rho B^{m-l}_{\mu'}f(x)|^2\om(x)\,dx
\\
\,\le\,&C_N \sum_{|\mu-\mu'|\le l}
\int| \Gamma_{\mu\mu'\rho}^{ml}
C^{m\delta}_\rho B^{m-l}_{\mu'}f(x)|^2w^{m,\delta}_{\mu\rho}*|\om|(x)
\,dx
\\
\,\le\,
&C_N \int|C^{m\delta}_\rho B^{m-l}_{\mu'}f(x)|^2
\sup_{|\mu-\mu'|\le l} w^{m,\delta}_{\mu\rho}*w^{m,\delta}_{\mu\rho}*|\om|(x)
\,dx.
\endalign$$
The  asserted inequality is an immediate consequence.\qed
\enddemo

We now conclude the proof of Proposition 4.6.  First,
since the 
maximal operator $M^{(1,2)}_m$ is invariant under two-parameter dilations
there is a scaled variant of Lemma 4.7. Also we can apply Lemma 4.7 twice,
in the $x_1-x_3$ and in the $x_2-x_3$ plane, the same applies to the 
scaled variant. We obtain
the inequality
$$
\sum_{k,\nu}\int
|V^{n\ep'}_{k\la}
A^{n}_{k\nu}f_k(x)|^2 
\om(x)\, dx\,\le\,
C \sum_{k,\nu'}\int 
|V^{n\ep'}_{k\la}
A^{n-l}_{k\nu'}f_k(x)|^2 \Cal M_n \om(x)\,dx
$$
if $l=(l_1,l_2)$ and $l_1\le n_1\ep'$, $l_2\le n_2\ep'$.
We iterate and  apply this inequality $N$ times; here $N\le 1+1/\ep'$.
The result is an estimate of the left hand side of (4.25) by an 
expression involving a scaled version of the square-function in Lemma 4.3
(with $A=\text{diag}(2^{k_1},2^{k_2})$). 
Namely if $\Gamma_{k\nu\la}^{n,\delta}$ is defined by
$$
\widehat{\Gamma_{k\nu\la}^{n,\delta}f}(\xi)\,=\,
\prod_{i=1}^2\bigl[
\chi^2(2^{m\delta}(2^{k_3}\xi_3-2^{-m\delta}\la_i))
\chi^2(2^{m\delta}(2^{k_i}\xi_i-2^{-m\delta}\nu_i))\bigr]\widehat f(\xi)
$$ 
we obtain the inequality
$$
\sum_{k,\nu}\int
|V^{n\ep'}_{k\la}
A^{n}_{k\nu}f_k(x)|^2 \om(x)\, dx\,\le\,
C_N \sum_{k,\nu,\la'}\int |\Gamma_{k\nu\la'}^{n,\ep'}
f_k(x)|^2\Cal M^{N}_n \om(x)\,dx,\qquad N\ge 1+\frac 1{\ep'},
$$
from which (4.25) follows by an application of Lemma 4.3.\qed
\enddemo
The asserted weighted norm inequality (4.10) now follows by an 
application of Lemma 4.5 and Proposition 4.6. This concludes the proof of 
Theorem 4.1.

\medskip

\remark{\bf Remark} The weighted inequality  in Proposition 3.6 implies
$$
\Bigl\|\Bigl(\sum_{k\nu}|A^n_{k\nu} \L_kf|^2\Bigr)^{1/2}\Bigr\|_4
\,\le\,C_\ep 2^{(n_1+n_2)\ep}\|f\|_4
$$
with $C_\ep=O(A^{1/\ep})$ as $\ep\to 0$, some $A>1$.
The geometrical arguments by C\'ordoba \cite{12} show that in fact
$C_\ep=O(\ep^{-a})$ for some $a>0$. It would be interesting 
 to find 
positive operators $\Cal N_\ep$, being uniformly bounded on $L^2$
 such that
$$
\sum_{k,\nu}\int
|A^{n}_{k\nu}\L_kf(x)|^2 \om(x)\, dx\,\le\,
C \ep^{-2a} 2^{(n_1+n_2)2\epsilon}
\int |f(x)|^2\Cal N_\ep[ \om](x)\,dx.
$$
An analogous problem is to find weighted norm inequalities for radial 
multipliers and associated maximal functions 
in $\Bbb R^2$, with a positive operator $\Cal N$.
In this context weighted  inequalities  with a nonpositive $\Cal N$
have  been proved in \cite{1}.
\endremark

\bigskip

\head
{\bf 5. $\boldkey{H}^{\boldkey p}$-estimates}
\endhead
The purpose of this section is to prove Theorem 1.7.
 The proof 
relies on a result on multiparameter Calder\'on-Zygmund theory 
obtained by the authors in
\cite{4} (extending earlier results by Journ\'e \cite{20} and  Fefferman 
\cite{17}). There it is shown for a large class of singular integral
operators $T$ that the boundedness of $T$ on certain scalar and
vector-valued rectangle atoms implies the boundedness on $H^p$.

To be precise let $R$ be an interval in $\Bbb R^d$
(\it i.e. \rm a rectangle parallel to the coordinate axes), and let $Q$ be 
a nonnegative integer. In what follows, $Q$ will always be $\ge
[1/p-1]$ (the largest integer $\le 1/p-1$).
Then $a$ is called a \it $(p,Q,R)$ rectangle atom \rm if 
$a$ 
is supported in $R$, if
$$
\int_R|a(x)|^2 dx\,\le\, |R|^{1-2/p}
$$
and if for $m=1,\dots,d$
$$\int_{\Bbb R^m}a(x_1,\dots,x_m,x_{m+1},\dots,x_d)x_1^{r_1}\dots x_m^{r_m}
\,dx_1\dots dx_m=0,\qquad r_1,\dots,r_m=0,\dots,Q$$
for almost  all $(x_{m+1},\dots,x_d)$; furthermore assume that the
analogous cancellation properties hold for all permutations of the variables 
$x_1,\dots,x_{d}$.

Now let $\Bbb R^d=\Bbb R^{d_1}\oplus\Bbb R^{d_2}$,
and let
$I$ be an interval in $\Bbb R^{d_1}$. Then we need the notion of
an $L^2(\Bbb R^{d_2})$-valued $(p,Q,I)$-rectangle 
atom. This is simply a function $a$ 
supported on $I\times\Bbb R^{d_2}$ such that
$$
\iint |a(x',x'')|^2 dx'' dx'\,\le\, |I|^{1-2/p}
$$
and
such that for $m=1,\dots,d_1$
$$\int_{\Bbb R^m}a(x_1,\dots,x_m,x_{m+1},\dots,x_{d_1+1},\dots,x_d)
 x_1^{r_1}\dots x_m^{r_m}
\,dx_1\dots dx_m=0, \qquad
r_1,\dots,r_m=0,\dots,Q$$
for almost  all $(x_{m+1},\dots,x_d)$; furthermore assume that the
analogous cancellation properties hold for all permutations of the variables 
$x_1,\dots,x_{d_1}$.

Now let $T:\,C_0^\infty(\Bbb R^d)\to (C_0^\infty(\Bbb R^d))'$ be an operator
with Schwartz kernel $K$, with the property that $K(x,y)$ is locally integrable
in $\{(x,y); x_i\neq y_i, i=1,\dots,d\}$.
Let $\Phi$ be a smooth bump function on $\Bbb R$ supported  in $[1,4]$
such that $\sum_{l=-\infty}^{\infty} \Phi(2^{-l}s)=1$ for $s\neq 0$.
For $\ell =(\ell_1,\dots,\ell_{d_1})$, $1\le d_1\le d$, define the operator
$T^\ell$ by
$$T^\ell f(x)=\int K(x,y) \prod_{i=1}^{d_1}
\Phi(2^{-\ell_i}|x_i-y_i|)\,f(y)\,dy.
$$

\proclaim{Theorem 5.1  \cite{4} }
Let  $0<p\le 1$, $s>d(d+1)/2$ and
 $Q\ge [1/p-1]$, $M>2$. Suppose that
\roster
\item $T$ is bounded on $L^2(\Bbb R^d)$ with operator norm $\le A$.
\item For all $d_1\in\{1,\dots,d-1\}$, for all $L\in \Bbb Z^{d_1}$, for all
intervals $I$ in $\Bbb R^{d_1}$ with sidelengths $2^{L_1},\dots,2^{L_{d_1}}$,
for all 
$L^2(\Bbb R^{d-d_1})$ valued 
$(p,Q,I)$ rectangle atoms $a$ and for all
$\ell=(\ell_1,\dots,\ell_{d_1})$,
$\ell_i>1$,
 $i=1,\dots,d_1$
$$\|T^{L+\ell}a\|_{L^p(\Bbb R^{d_1},
L^2(\Bbb R^{d_2}))}\,\le\,A(\sum_{i=1}^{d_1}\ell_i)^{-s/p}.
\tag 5.1
$$
\item The condition analogous to {\rm (5.1)} is valid for every permutation of 
the variables $x_1,\dots,x_d$.
\item For all  $L\in \Bbb Z^{d}$, for all
intervals $R$ in $\Bbb R^{d}$ with sidelengths $2^{L_1},\dots,2^{L_{d}}$,
for all $(p,Q,R)$ rectangle atoms $a$ and for all 
$\ell=(\ell_1,\dots,\ell_{d})$,
$\ell_i>1$, $i=1,\dots,d-1$
$$\|T^{L+\ell}a\|_{L^p(\Bbb R^{d})}\,\le\,A
(\sum_{i=1}^d\ell_i)^{-s/p}
\tag 5.2
$$
\endroster
Then $T$ extends to a  bounded operator
 from the multiparameter Hardy space $H^p(\Bbb R^d)$  to 
$L^p(\Bbb R^d)$ and the operator norm is bounded by $C\,A$. Here $C$ depends 
only on
$p$, $d$ and  $s$.
If $T$ 
is  translation invariant  then
$T$ is bounded on $H^p(\Bbb R^d)$.
\endproclaim

We now consider convolution operators $T$
given by Fourier multipliers $m$ via
$\widehat{Tf}(\xi)=m(\xi)\widehat f(\xi)$.
For $k\in\Bbb Z^{d_1}$ let $T_k$ be the operator with Fourier multiplier
$m(\xi)\prod_{i=1}^{d_1} \beta(2^{k_i}\xi_i)$. 
Variants of the standard Marcinkiewicz multiplier theorem on $H^p$ spaces 
follow from  Theorem 5.1 and

\proclaim{Proposition 5.2}
Suppose that $0<p\le 1$, $\alpha>1/p-1/2$
and let $Q$, $\ep$ be such that
$Q\ge[1/p-1]$ and $0<2\epsilon<\min\{\alpha-1/p+1/2,Q-1/p+2,1\}$.
\roster
\item
Suppose that $1\le d_1\le d-1$ and
$$\sup_{t\in(\Bbb R_+)^{d_1}}\sup_{(\xi_{d_1+1},\dots,\xi_d)\in 
\Bbb R^{d-d_1}}
\|\beta_{(1)}\otimes\dots\otimes\beta_{(d_1)}
\,m(t_1\cdot,\dots,t_{d_1}\cdot,\xi_{d_i+1},
\dots,\xi_d)\|_{\Cal H^2_\alpha(\Bbb R^{d_1})}
<\infty.
\tag 5.3
$$
Then  for all $L\in \Bbb Z^{d_1}$, for all
intervals $I$ in $\Bbb R^{d_1}$ with sidelengths $2^{L_1},\dots,2^{L_{d_1}}$,
for all 
$L^2(\Bbb R^{d-d_1})$ valued 
$(p,Q,I)$ rectangle atoms $a$, for all
$\ell=(\ell_1,\dots,\ell_{d_1})$,
$\ell_i>1$, $i=1,\dots,d_1$ and for all $k\in \Bbb Z^{d_1}$
$$\|(T_k)^{L+\ell}a\|_{L^p(\Bbb R^{d_1},
L^2(\Bbb R^{d-d_1}))}\,\le\,C\,A\prod_{i=1}^{d_1}
2^{-\epsilon (\ell_i+|k_i|)}.
\tag 5.4
$$
\item 
The inequality analogous to {\rm (5.4)} holds for every permutation of 
the variables $x_1,\dots,x_d$.
\item
Suppose that
$$\sup_{t\in(\Bbb R_+)^d}
\|\beta_{(1)}\otimes\dots\otimes\beta_{(d)}
\,m(t_1\cdot,\dots,t_d\cdot)
\|_{\Cal H^2_\alpha(\Bbb R^{d})}
<\infty.
\tag 5.5
$$
Then for all  $L\in \Bbb Z^{d}$, for all
intervals $R$ in $\Bbb R^{d}$ with sidelengths $2^{L_1},\dots,2^{L_{d}}$,
for all $(p,Q,R)$ rectangle atoms $a$, for all 
$\ell=(\ell_1,\dots,\ell_{d})$,
$\ell_i>1$, $i=1,\dots,d$, for all $k\in\Bbb Z^d$
$$\|(T_k)^{L+\ell}a\|_{L^p(\Bbb R^{d})}\,\le\,C\,A\prod_{i=1}^{d}
2^{-\epsilon (\ell_i+|k_i|)}.
\tag 5.6
$$
\endroster
If {\rm(5.5)} is valid then $m$ is bounded and  {\rm(5.3)} and
the analogous conditions obtained by permuting variables are also satisfied.
In particular {\rm(5.5)} implies that
 $T$ is bounded on the multiparameter Hardy space 
$H^p(\Bbb R^d)$   and the operator norm is bounded by $C\,A$. 
\endproclaim

Proposition 5.2 is proved by standard arguments, see for example 
the proof of \cite{4, Proposition 5.1}.
The last conclusion of Proposition 5.2 follows of course  by Theorem 5.1.
The reader should note that the multipliers in Theorem 1.7 generally do not
satisfy the assumption (5.5), even in the two-dimensional case.

\demo{\bf Proof of Theorem 1.7}
We may clearly assume that $p\le 1$.
Again since characteristic 
functions of  half spaces with boundaries  parallel to 
the coordinate axes are Fourier multipliers of multiparameter Hardy spaces
there is no loss of generality in assuming that $m$ is supported in
$\{\xi;\,\xi_i\ge 0,\,i=1,\dots,d\}$. 
We use the notations introduced in the proof of Theorem 4.1.
Let $T^n_k$ is as in (4.9) and set
 $T^n=\sum_{k\in\Bbb Z^d}T_k^n$. We shall show that
$T^n$ is bounded on $H^p(\Bbb R^d)$ with operator norm bounded by
$$
C_s \sup_{k\in\Bbb Z^d}\|g^n_k\|_{L^p(\Bbb R^{d-1})} 2^{(n_1+\dots+n_{d-1})
(\frac 2p-1)}
(1+n_1+\dots+n_{d-1})^{(s+d)/p}.
\tag 5.7
$$
Since
$$\multline
\sum_{n\in(\Bbb N_0)^{d-1}}
\sup_{k\in\Bbb Z^d}\|g^n_k\|_{L^p(\Bbb R^{d-1})} 
2^{(n_1+\dots+n_{d-1})(\frac 2p-1)}
(1+n_1+\dots+n_{d-1})^{(s+d)/p}
\\ \,\,\,\,\le\,
C_\ep\sup_{k\in\Bbb Z^d}\|g_k\|_{\Cal H^p_\alpha}
\qquad\text{if }\alpha>2/p-1
\endmultline
$$
the conclusion of  Theorem 1.7 follows.

We have to verify  the hypotheses (5.1) and (5.2) of 
Theorem 5.1 for the operator $T^n$. 
The mixed norm inequalities are a straightforward  consequence of Proposition 
5.2. In order to see this let
$$F_h(\xi)=h(\xi_1/\xi_d,\dots,\xi_{d-1}/\xi_d)$$
where $h$ is compactly supported in $[1/2,2]^{d-1}$.
Then for $\alpha\ge 0$ one has the inequalities
$$
\sup_{\xi_d}
\|\beta_{(1)}\otimes\dots\otimes\beta_{(d)}\,F_h(\cdot,\xi_d)
\|_{\Cal H^2_\alpha
(\Bbb R^{d-1})}
\,\le\,C\|h\|_{\Cal H^2_\alpha(\Bbb R^{d-1})}.
\tag 5.8
$$
and
$$
\multline
\sup_{\xi_1}
\|\beta_{(1)}\otimes\dots\otimes\beta_{(d)}\,F_h(\xi_1,\cdot)
\|_{\Cal H^2_\alpha(\Bbb R^{d-1})}
\\
\,\le\,C
\cases
\|h\|_{\Cal H^2_\alpha(\Bbb R)}\qquad&\text{if }d=2
\\ 
\|h\|_{\Cal H^2_\alpha(\Bbb R^{d-1})}
+\sup_{s_1}\sum_{k=2}^{d-1}
\|\Cal D_2^\alpha\dots\Cal D_{k-1}^\alpha\Cal D_k^{2\alpha}\Cal D_{k+1}^\alpha
\dots\Cal D^\alpha_{d-1} h(s_1,\cdot)\|_{L^2(\Bbb R^{d-2})}
\qquad&\text{if }d\ge 3.
\endcases
\endmultline
\tag 5.9$$
It is straightforward to verify (5.8) and (5.9) if $\alpha$ is 
a nonnegative even integer and the general case
follows by analytic interpolation. 
Note also that by a version of
Sobolev's imbedding theorem
$\Cal H^p_\beta(\Bbb R^{d-1})\subset
\Cal H^2_\alpha(\Bbb R^{d-1})$ if $p\le 2$ and
$\beta\ge \alpha+(1/p-1/2)$. Using this and (5.8), (5.9)
we see that (5.3) is verified for the case $d_1=d-1$. The other cases 
follow similarly. An application of Proposition 5.2 implies (5.2).

The main work in the proof  consists in the verification of (5.1). 
Assume that
$a$ is a $(p,Q, R) $ rectangle atom and $R$ is an interval of dimensions
$2^{L_1}\times \dots\times 2^{L_{d}}$.
Then we
shall prove that
$$
\|(T_k^n)^{L+\ell}a\|_p\,\le\,C 
2^{(n_1+\dots+n_{d-1})N} \prod_{i=1}^d 2^{-\ep(\ell_i+|k_i-L_i|)}\,
\|g_k^n\|_\infty,\qquad N>2(1/p-1/2)
\tag 5.10
$$
for some $\ep>0$ and also
$$
\|(T^n_k)^{L+\ell}a\|_p\,\le\,C 
2^{(n_1+\dots+n_{d-1})(2/p-1)}\|g_k^n\|_p.
\tag 5.11
$$
We shall use (5.11) only if $\max_j\{k_j-L_j\},\max_j\{\ell_j\}\le C_p
(1+\sum_i n_i)$ where $C_p$ is a large fixed constant while
(5.10) is a remainder estimate.
 In fact
applying the Sobolev inequality (4.8) with $d_1=0$  we see that
(5.10) and (5.11) imply
$$
\align
&\Bigl\|\sum_{k\in \Bbb Z^d}(T^n_k)^{L+\ell}a\Bigr\|_p\,
\\
&{\aligned\le\,C\biggl(
&\sum\Sb\max\{|k_i-L_i|,i=1,\dots,d\}\ge
\\\ep^{-1}(2N+2/p)(n_1+\dots+n_{d-1})\endSb
2^{(n_1+\dots+n_{d-1})(Np+1)} \prod_{i=1}^d 2^{-\ep p(\ell_i+|k_i-L_i|)}\,
\|g^n_k\|_p^p
\\+&\sum\Sb\max\{|k_i-L_i|,i=1,\dots,d\}
<\\ \ep^{-1}(2N+2/p)(n_1+\dots+n_{d-1})\endSb
\min\bigl\{2^{(n_1+\dots+n_{d-1})(2-p)};
2^{(n_1+\dots+n_{d-1})(Np+1)} \prod_{i=1}^d 2^{-\ep p(\ell_i+|k_i-L_i|)}
\bigr\}
\|g^n_k\|_p^p\biggr)^{1/p}\endaligned}
\\
\,&\le\, C
\,2^{(n_1+\dots+n_{d-1})(2/p-1)}
\frac{(1+n_1+\dots+n_{d-1})^{(s+d)/p}}
{(\ell_1+\dots+\ell_d)^{s/p}}
\|g_k^n\|_p
\endalign
$$
and it follows that $T^n$ is bounded on $H^p$ with norm not exceeding (5.7).

The verification of (5.10) is  easy.  Simply observe that
$$
\bigl|\partial_\xi^\gamma [\phi(\xi)m_{k,n}(\xi)]\bigr|\,\le\,C_\gamma
2^{n_1\gamma_1+\dots+
n_{d-1}\gamma_{d-1}} (2^{n_1\gamma_d}+\dots+2^{n_{d-1}\gamma_d})
$$
and an application of Proposition 5.2 yields (5.10).

We now verify (5.11) and assume for convenience $d=3$. 
We show that
 (using the notation introduced in (4.15)) 
$$
\align
\|(T^n_{k\nu})^{L+\ell}a\|_p\,&\le\,C 
b^n_{k\nu}\,2^{(n_1+n_2)(1/p-1)}
\|g_k^n\|_p
\tag 5.12
\\
\|(T^{n,1}_{k\nu_1})^{L+\ell}a\|_p\,&\le\,C_N
b^{n,1}_{k\nu_1}\,2^{n_1(1/p-1)}2^{-n_2 N}
\|g_k^n\|_p
\tag 5.13
\\
\|(T^{n,2}_{k\nu_2})^{L+\ell}a\|_p\,&\le\,C_N
b^{n,2}_{k\nu_2}\,2^{-n_1 N}2^{n_2(1/p-1)}
\|g_k^n\|_p
\tag 5.14
\\
\|(T^{n,0}_k)^{L+\ell}a\|_p\,&\le\,C_N
b^n_{k}\,2^{-(n_1+n_2)N}
\|g_k^n\|_p
\tag 5.15
\endalign
$$
Using (4.16) and (4.17) with $r=p$ we see that (5.11) follows from (5.12-15).
We shall only verify (5.12); the remaining cases are similar or simpler.

We divide the rectangle $R$ (which has dimensions 
$2^{L_1}\times 2^{L_2}\times 2^{L_3}$) into
$\prod_{i=1}^3\max\{1,2^{L_i-k_i}\}$  congruent
intervals $R_k^\mu$ of dimensions
$$\min\{2^{L_1},2^{k_1}\}\times \min\{2^{L_2},2^{k_2}\}\times
\min\{2^{L_3},2^{k_3}\}
$$
and centers $y^\mu_k$.
Let be $a_{k}^\mu=a\chi_{R_k^\mu}$ and
let
$$
\Cal R^{L+\ell}_\mu\,=\,\{x;\,2^{L_i+\ell_i-2}\le|x_i-(y^\mu_k)_i|\le 2^{L_i+
\ell_i+2},
i=1,2,3\}.
$$
Then it is easy to check that if 
$y\in \supp a_k^\mu$, $x\in\Cal R^{L+\ell}_\mu$
then for $U_{k\nu, N}^n$ as in Lemma 4.4
$$
U_{k\nu,N}^n(x-y)\approx
U_{k\nu,N}^n(x-y_k^\mu )
$$
and therefore by Lemma 4.4
$$
\align
&\|(K^n_{k\nu}\Phi_{L+\ell})*a\|_p
\le\,C\,
\bigl(\sum_\mu\|(K^n_{k\nu}\Phi_{L+\ell})*a_{k}^\mu\|_p^p\bigr)^{1/p}
\\
&\le\,C_N\,b^n_{k\nu}
\Bigl(\sum_\mu\int
\Bigl[\int U^n_{k\nu,N}(x-y)
|\Phi_{L+\ell}(x-y)|
|a_{k}^{\mu}(y)|\,dy\,\Bigr]^p\,dx\Bigr)^{1/p}
\\
&{\aligned\le\,
C_N\,2^{(n_1+n_2)(\frac 1p-1)}\,b^n_{k\nu}
\Bigl(\sum_\mu
\int_{\Cal R^{L+\ell}_\mu} 
&2^{-(n_1+n_2+k_1+k_2+k_3)(p-1)}\bigl[U^n_{k\nu N}(x-y_k^\mu )\bigr]^{p}dx 
\\
\,&\times\Bigl[\int |a_k^\mu(y)|dy\,
2^{(k_1+k_2+k_3)(\frac 1p-1)}\Bigr]^p\Bigr)^{1/p}.
\endaligned}
\tag  5.16
\endalign
$$
Using H\"older's inequality we see that
$$\multline
2^{(L_1+L_2+L_3)(\frac 1p-1)}
\Bigl(
\sum_\mu 
\prod_{i=1}^3\bigl[\min\{1,2^{(k_i-L_i)(1-p)}\}\bigr]
\|a_k^\mu\|_1^p\Bigr)^{1/p}\\
\,\le\,C\,2^{(L_1+L_2+L_3)(\frac 1p-1)}\|a\|_1\,\le\,C\,
|R|^{\frac 1p-\frac 12}\|a\|_2\,\le\,C'.
\endmultline
\tag 5.17
$$
We
perform the linear volume preserving  change of variables
$$
(v_1,v_2,v_3)\,=\,(x_1,x_2,2^{k_1-k_3}u^1_{\nu_1}x_1
+2^{k_2-k_3}u^2_{\nu_2}x_2+x_3)
$$
and see that for $N>1/p$
$$
\multline
\int_{\Cal R^{L+\ell}_\mu} 2^{-(n_1+n_2+k_1+k_2+k_3)(p-1)}
\bigr[U^n_{k\nu, N}(x-y_k^\mu)\bigl]^{p}\,dx 
\\
\,\le\,C\int
\frac{2^{-k_1-n_1}}{(1+|2^{-k_1-n_1} v_1|)^{Np}}
\frac{2^{-k_2-n_2}}{(1+|2^{-k_2-n_2} v_2|)^{Np}}
\frac{2^{-k_3}}{(1+|2^{-k_3} v_3|)^{Np}}dv
\,\le\,C'.
\endmultline
$$
Therefore if 
$k_i-L_i\le 0$, $i=1,2,3$, the desired estimate (5.12) follows from (5.16) 
and (5.17).

In all other cases we use similar arguments together
with the  cancellation properties of the atom.  
For example assume $k_1\le L_1$, $k_2\le L_2$, $k_3\le L_3$.
Since
$$
\iint a^\mu_k(y_1,y_2,y_3)y_1^{r_1}y^{r_2}_2\,dy_1dy_2\,=\,0
$$
for almost all $y_3$ for $0\le r_1,r_2\le Q$ we see using Taylor's formula
that
$$
\multline
(K^n_{k\nu}\Phi_{L+\ell})*a_k^\mu(x_1,x_2,x_3)
\,=\,
\\
\int_0^1\frac{(1-s)^Q}{Q!}
\int\Bigl(\frac{\partial}{\partial x_3}\Bigr)^{Q+1}(K^n_{k\nu}\Phi_{L+\ell})
(x_1-y_1,x_2-y_2,x_3-(y_k^\mu)_3+s((y_k^\mu)_3-y_3))
((y_k^\mu)_3-y_3)^{Q+1}
a^\mu_k(y)\,dy\, ds
\endmultline
$$
and using Leibniz' rule and Lemma 4.4 we see that
$$
\big|(K^n_{k\nu}\Phi_{L+\ell})*a_k^\mu(x_1,x_2,x_3)\big|
\,\le\,C
2^{L_3(Q+1)}\max\{2^{-k_3(Q+1)}, 2^{-(L_3+\ell_3)(Q+1)}\}
U^n_{k\nu, N}(x-y_k^\mu)
\, b^n_{k\nu}\,\|a^\mu_k\|_1.
$$

Similar considerations in the other cases (where we  use that
$a^\mu_k$ has cancellation in the $y_i$ variable whenever
$k_i\ge L_i$)  lead to
$$
\multline
\|(K^n_{k\nu}\Phi_{L+\ell})*a\|_p
\,
\le\,C_N\,
\prod_{i=1}^3\bigl[\min\{1,(2^{-\ell_i}+2^{L_i-k_i})\}\bigr]^{Q+1}
\\
\times\,b^n_{k\nu}
\Bigl(\sum_\mu
\int_{\Cal R^{L+\ell}_\mu} 
2^{(n_1+n_2+k_1+k_2+k_3)(p-1)}
\bigr[U^n_{k\nu, N}(x-y_k^\mu)\bigl]^{p}dx 
\|a_k^\mu\|_1^p
\Bigr)^{1/p}.
\endmultline
$$
As above it is easy to check that for $N>1/p$
$$
\multline
\int_{\Cal R^{L+\ell}_\mu} 2^{-(n_1+n_2+k_1+k_2+k_3)(p-1)}
\bigr[U^n_{k\nu, N}(x-y_k^\mu)\bigl]^{p}\,dx 
\\
\le\,C\,
\min\{1,2^{L_1+\ell_1-k_1-n_1}\}
\min\{1,2^{L_2+\ell_2-k_2-n_2}\}
\min\{1,2^{L_3+\ell_3-k_3}\}.
\endmultline
$$
Therefore
$$
\alignat2
\|(K^n_{k\nu}\Phi_{L+\ell})*a\|_p
\,\le\,
&C2^{(n_1+n_2)(\frac 1p-1)}&&b^n_{k\nu}
\prod_{i=1}^3\bigl[
\min\{1,2^{(L_i+\ell_i-k_i)/p}\}
\min\{1,(2^{-\ell_i}+2^{L_i-k_i})^{Q+1}\}
\bigr]
\\
&{} &&\times 2^{(k_1+k_2+k_3)(\frac 1p-1)}
\Bigl(\sum_\mu\|a^\mu_k\|_1^p
\Bigr)^{1/p}
\\
\,\le\,&
C2^{(n_1+n_2)(\frac 1p-1)}&&
b^n_{k\nu}
2^{(L_1+L_2+L_3)(\frac 1p-1)}
\Bigl(\sum_\mu\|a^\mu_k\|_1^p
\prod_{i=1}^3\bigl[\min\{1,2^{(L_i-k_i)(\frac 1p-1)}\}\bigr]^p
\Bigr)^{1/p}
\\
\,\le\,&C2^{(n_1+n_2)(\frac 1p-1)}&&b^n_{k\nu}
|R|^{\frac 1p-1}\|a\|_1
\\
\,\le\,&C2^{(n_1+n_2)(\frac 1p-1)}&&b^n_{k\nu}.
\endalignat
$$
This proves (5.12) and concludes the proof of Theorem 1.7.\qed
\enddemo

\end
\bigskip



\Refs\nofrills{{\bf References}}

\ref \no 1 \by A. Carbery \paper A weighted inequality for the maximal 
Bochner-Riesz operator on $R^2$
\jour Trans. Amer. Math. Soc. \vol 287 \yr 1985 \pages 673--679
\endref

\ref \no 2 \bysame
\paper Differentiation in lacunary directions and an extension of the 
Marcin\-kiewicz multiplier theorem
\jour Ann. Inst. Fourier
(Grenoble)
\vol 38\yr 1988\pages 157--168\endref



\ref \no 3 \by A. Carbery and A. Seeger
 \paper Conditionally convergent
series of linear operators in $L^p$ spaces and $L^p$ estimates for
pseudo-differential operators \jour Proc. London Math. Soc. \vol 57
\yr 1988 \pages 481--510
\endref


\ref\no  4 
\bysame
  \paper $H^p$ and $L^p$ variants of
 multiparameter Calder\'on-Zygmund theory \jour Trans. Amer. Math. Soc.
\vol 334\yr 1992\pages 719--747
\endref 


\ref\no 5\bysame\paper 
Multiparameter Calder\'on-Zygmund theory and
interpolation of analytic families 
\jour preprint\yr 1993
\endref

\ref \no 6 \by S.Y.A. Chang and R. Fefferman
\paper The Calder\'on-Zygmund decomposition on product domains
\jour Amer. J. Math \vol 104 \yr 1982 \pages 445--468 \endref

\ref\no 7\by M. Christ\paper On three related problems in harmonic analysis
\jour Unpublished manuscript\yr 1984\endref

\ref \no 8 \by R. R. Coifman, J. L. Rubio de Francia and S. Semmes
\paper Multiplicateurs de Fourier dans $L^p(\Bbb R)$ et estimations
quadratiques \jour C. R. Acad. Sci. Paris
\vol 306 \yr 1988 \pages 351--354
\endref

\ref\no 9\by W. C. Connett and A. L. Schwartz\paper A remark about Calder\'on's
upper $s$ method of interpolation
\inbook Interpolation spaces and allied topics in analysis
\bookinfo Proceedings, Lund  1983, ed. by M. Cwikel and J. Peetre
\pages 48--53
\publ
Lecture notes in Math. 1070, Springer-Verlag
\publaddr Berlin, Heidelberg \yr 1984
\endref


\ref\no 10 \by A. C\'ordoba
\paper A note on Bochner-Riesz operators \jour Duke Math. J.
\vol 46 \yr 1979 \pages 505-511
\endref

\ref\no 11 \bysame
\paper Maximal functions, covering lemmas and Fourier multipliers
\inbook Harmonic analysis in Euclidean spaces
\bookinfo Proceedings of symposia in pure mathematics,
vol. XXXV, part 1, ed. by S. Wainger and G. Weiss\pages 29--50
\publ American Mathematical Society\publaddr Providence\yr
1979\endref

\ref\no 12 \bysame
\paper  Geometric Fourier analysis 
\jour Ann. Inst. Fourier
\vol  32 \yr 1982
\pages 215--226
\endref

\ref\no 13\by  A. C\'ordoba and C. Fefferman\paper A weighted norm inequality
for singular integrals\jour Studia Math.\vol 57\yr 1976\pages 97--101
\endref

\ref\no 14\by  A. C\'ordoba and R. Fefferman\paper On the equivalence
between the boundedness of certain classes of maximal and multiplier 
operators in Fourier analysis
\jour Proc. Nat. Acad. Sci. USA
\vol 74 \yr 1977\pages 423--425
\endref

\ref\no 15 \by J. Duoandikoetxea and A. Moyua
\paper Homogeneous multipliers in the plane
\jour Proc. Amer. Math. Soc.
\vol 112\yr 1991\pages 441--450
\endref

\ref\no 16\by C. Fefferman\paper The multiplier problem for the ball
\jour Annals of Math.\vol 94\yr 1971\pages 331--336
\endref

\ref \no 17\by R. Fefferman \paper Harmonic analysis on
product spaces \jour Annals of Math. \vol 126 \yr 1987
\pages 109--130
\endref



\ref \no 18
\by J. Garc\'\i a-Cuerva and J. L. Rubio de Francia
\book Weighted norm inequalities and related topics
\bookinfo North-Holland Math. Studies {\bf 116}
\publ North-Holland\publaddr Amsterdam, New York, Oxford
\yr 1985
\endref

\ref \no 19 \by C. Herz and N. M. Rivi\`ere
\paper Estimates for translation invariant
operators on spaces with mixed norm
\jour Studia Math.
\vol 44 \yr 1972 \pages 511--515
\endref

\ref \no 20
\by J.-L. Journ\'e\paper Calder\'on-Zygmund operators on product spaces
\jour Rev. Mat. Iberoamericana \vol 1 \yr 1985 \pages 55--91
\endref

\ref\no 21\by C. Kenig and P. Tomas\paper $L^p$ behavior of certain second 
order partial differential operators
\jour Trans. Amer. Math. Soc.\vol 262\yr 1980\pages 521--531\endref


\ref \no 22\by B. L\'opez-Melero
\book An\'alisis arm\'onico en $\Bbb R^n$: tres escenarios
\bookinfo Tesis Doctoral
\publaddr Universidad Aut\'onoma, Madrid \yr 1981
\endref


\ref \no  23\by A. Nagel, E. M. Stein and S. Wainger
\paper Differentiation in lacunary directions
\jour Proc. Nat. Acad. Sci. USA \vol 75 \yr 1978 \pages 1060--1062
\endref


\ref\no 24\by J.-L. Rubio de Francia\paper Estimates 
for some square functions of
Littlewood-Paley type\jour Publicacions  Mathem\`atiques\vol 27\yr 1983\pages 
81--108\endref


\ref\no  25\by E. M. Stein   \book Singular integrals and differentiability
properties of functions \publ Princeton Univ. Press \publaddr Princeton, N.J.
\yr 1971 \endref


\ref\no 26\by H. Triebel\paper Eine Bemerkung zur 
nicht-kommutativen Interpolation\jour Math. Nachr.\vol 69\pages 57--60\yr 1975
\endref

\ref\no  27 \bysame
 \book Theory of function spaces 
\publ Birkh\"auser Verlag \publaddr Basel, Boston, Stuttgart
\yr 1983
\endref



\endRefs
